\documentclass{article}

\usepackage{amsmath}
\usepackage{amsthm}
\usepackage{amssymb}
\usepackage{mathrsfs}
\usepackage{graphicx,color}
\usepackage{texdraw}

\newtheoremstyle{plain2}{\topsep}{\topsep}%
     {\itshape}
     {}
     {\bfseries}
     {.}
     {.5em}
     {\thmnumber{(#2)}\thmname{ #1}\thmnote{ #3}}

\theoremstyle{plain2}
\newtheorem{teo}{Theorem}[section]
\newtheorem{prop}[teo]{Proposition}
\newtheorem{coro}[teo]{Corollary}
\newtheorem{lemma}[teo]{Lemma}

\newtheoremstyle{definition2}{\topsep}{\topsep}%
     {}
     {}
     {\bfseries}
     {.}
     {.5em}
     {\thmnumber{(#2)}\thmname{ #1}\thmnote{ #3}}

\theoremstyle{definition2}

\newtheorem{rem}[teo]{Remark}
\newtheorem{defi}[teo]{Definition}

\newcommand{\parder}[2]{\frac{\partial #1}{\partial #2}}

\def\C{\mathbb{C}}
\def\R{\mathbb{R}}
\def\Z{\mathbb{Z}}
\def\a{\alpha}
\def\theta{\vartheta}
\def\phi{\varphi}
\newcommand{\curl}{\nabla\times}
\renewcommand{\div}{\nabla\cdot}

\newcommand{\bfdx}{\textrm{d}\mathbf{x}}
\newcommand{\xbar}{\bar{x}}

\newcommand{\mpA}{\mathbf{A}}
\newcommand{\mfB}{\mathbf{B}}

\newcommand{\Sgamma}{\mathcal{S}_\gamma}
\newcommand{\ui}{u_i}
\newcommand{\uj}{u_j}
\newcommand{\hatui}{\hat{u}_i}
\newcommand{\HAD}{\mathcal{H}_{\mpAa}^1(D)}

\newcommand{\Da}{D_a}
\newcommand{\mpAa}{\mathbf{A}_a}
\newcommand{\Omegaa}{\Omega_a}
\newcommand{\Sigmaa}{\Sigma_a}
\newcommand{\Ua}{U_a}
\newcommand{\UO}{U_0}
\newcommand{\Thetaa}{\Theta_a}
\newcommand{\ThetaO}{\Theta_0}
\newcommand{\mpAO}{\mathbf{A}_0}
\newcommand{\thetaa}{\theta_a}
\newcommand{\Gammaa}{\Gamma_a}
\newcommand{\GammaO}{\Gamma_0}

\newcommand{\Ga}{G_{1,a}}
\newcommand{\Fa}{F_a}

\newcommand{\etaa}{\eta_a}

\newcommand{\rhoa}{\rho_a}
\newcommand{\rhoO}{\rho_0}
\newcommand{\phia}{\phi_a}
\newcommand{\phiO}{\phi_0}

\newcommand{\etaunoa}{\eta_a^{(1)}}
\newcommand{\etaunoO}{\eta_0^{(1)}}
\newcommand{\xitildea}{\tilde{\xi}_a}
\newcommand{\xia}{\xi_a}
\newcommand{\xiO}{\xi_0}

\newcommand{\ua}{u_a}

\newcommand{\uO}{u_0}

\newcommand{\Omegatilde}{\tilde{\Omega}}
\newcommand{\sigmatilde}{\tilde{\sigma}}
\newcommand{\mpAtilde}{\tilde{\mathbf{A}}}
\newcommand{\Gammatilde}{\tilde{\Gamma}}

\newcommand{\uuno}{u^{(1)}}
\newcommand{\uunoa}{u^{(1)}_a}
\newcommand{\uunoO}{u^{(1)}_0}
\newcommand{\gammauno}{\gamma^{(1)}}

\newcommand{\Vuno}{V^{(1)}}
\newcommand{\Vunoa}{V^{(1)}_a}
\newcommand{\VunoO}{V^{(1)}_0}

\newcommand{\udue}{u^{(2)}}
\newcommand{\uduea}{u^{(2)}_a}
\newcommand{\udueO}{u^{(2)}_0}

\newcommand{\gammadue}{\gamma^{(2)}}
\newcommand{\gammaduea}{\gamma^{(2)}_a}
\newcommand{\gammadueO}{\gamma^{(2)}_0}
\newcommand{\Vdue}{V^{(2)}}
\newcommand{\Vduea}{V^{(2)}_a}
\newcommand{\VdueO}{V^{(2)}_0}

\newcommand{\za}{z_a}
\newcommand{\Za}{Z_a}
\newcommand{\uh}{u^{(h)}}
\newcommand{\wh}{w^{(h)}}
\newcommand{\uinfty}{u^\infty}
\newcommand{\winfty}{w^\infty}
\newcommand{\Gzero}{\mathcal{G}_{0}}

\title{\sc Nodal sets of magnetic Schr\"odinger operators of Aharonov--Bohm type and energy minimizing partitions
\footnote{Work partially supported by MIUR, Project ``Metodi
Variazionali ed Equazioni Differenziali Non Lineari'' }
}

\author{
Benedetta Noris and Susanna Terracini \\
}

\begin{document}

\maketitle

\begin{abstract}
In this paper we consider a stationary Schr\"odinger operator in the plane, in presence of a magnetic field of Aharonov--Bohm type with semi--integer circulation. We analyze the nodal regions for a class of solutions such that the nodal set consists of regular arcs, connecting the singular points with the boundary. In case of one magnetic pole, which is free to move in the domain, the nodal lines may cluster dissecting the domain in three parts. Our main result states that the magnetic energy is critical (with respect to the magnetic pole) if and only if such a configuration occurs. Moreover the nodal regions form a minimal 3--partition of the domain (with respect to the real energy associated to the equation), the configuration is unique and depends continuously on the data. The analysis performed is related to the notion of spectral minimal partition introduced in \cite{HHOT}. As it concerns eigenfunctions, we similarly show that critical points of the Rayleigh quotient correspond to multiple clustering of the nodal lines.
\bigskip


\end{abstract}

\section{Introduction}
According to the recent literature, the analysis of Aharonov--Bohm operators with half integer circulation may lead to new insights on the nodal configuration of eigenfunctions (see \cite{BHH}) to the Dirichlet laplacian in planar domains. A useful related concept is that of \emph{spectral minimal partition}, introduced in \cite{HHOT, HHOT2} (see also \cite{H} for a survey and \cite{HHOT3} for further developments), that is a partition which is optimal with respect to the largest of the first eigenvalues. According to the analysis there, spectral minimal partitions need not to be nodal (i.e. nodal partitions associated to an eigenfunction), for the components may cluster in an odd number. In the planar case, a natural way to handle triple clustering of subdomains is to pass on the double covering of the punctured domain. Equivalently, we can associate a new operator of Aharonov--Bohm type by introducing a singular magnetic field with pole at the triple (or multiple) nodal junction, having care of prescribing half integer circulation (see \cite{BHH,HHO,HHO1}). In particular this shows that selected solutions of the associated Aharonov--Bohm equation, though complex, do possess nodal sets having interesting spectral properties. 

In this paper we reverse this point of view and, starting from a solution to an Aharonov-Bohm operator, we study its nodal lines. Next we move the pole of the magnetic potential and study the dependence of the nodal lines with respect to such a singularity. We are specially interested in multiple clustering. Our main result relates the  occurrence of triple nodal junctions  with criticality of the corresponding magnetic energy. A strictly related statement is the following.
\begin{teo}\label{theorem_eigenvalues}
Let $\lambda_a^k$ ($k$ positive integer) be the $k$--eigenvalue of the Aharonov--Bohm operator with half integer circulation and Dirichlet boundary condition in a planar domain $\Omega$. Here $a$ is the singularity of the magnetic potential. Assume that the function $a\mapsto \lambda_a^k$ is differentiable and that its gradient vanishes at $a\in\Omega$. Then the corresponding eigenfunctions possess an odd number, greater than or equal to three, of nodal arcs intersecting at $a$.
\end{teo}

As a consequence, we have the following:

\begin{coro}\label{coro_eigenvalues} If $a\in\Omega$ is a local extremum point of a simple eigenvalue $\lambda^k_a$, then there is a multiple junction of at least three nodal lines of the corresponding eigenfunctions at $a$. The first eigenvalue function $a\mapsto \lambda^1_a$ has a global interior maximum (where it is not differentiable) corresponding to an eigenfunction of multiplicity exactly two.
\end{coro}

This results complement the upper bound on the multiplicity of the first eigenvalues given in \cite{HHO1} and provides a theoretical key for the explanation of many of the numerical results of \cite{BHH,BHV,BH}.

Let us state more precisely the setting and results of the paper. Let $\Omega\subset \R^2$ be a simply connected domain with regular boundary. Given a point $a \in \Omega$ we consider the stationary magnetic Schr\"odinger operator
\begin{equation*}
H_{\mpAa,V}=(i \nabla + \mpAa)^2 + V,
\end{equation*}
acting on complex valued functions $U \in L^2(\Omega, \C)$. Here $V(x)\in W^{1,\infty}(\Omega)$ represents a conservative potential, whereas the magnetic potential has the form
\begin{equation}\label{potential}
\mpAa(x_1,x_2)=\frac{2n+1}{2} \left( - \frac{x_2-a_2}{(x_1-a_1)^2+(x_2-a_2)^2}, \frac{x_1-a_1}{(x_1-a_1)^2+(x_2-a_2)^2} \right) +\nabla\Phi
\end{equation}
where $a=(a_1,a_2), n\in \Z$ and $\Phi$ is any function\footnote{The vector potential is determined up to gauge transformations.} of class $C^2(\Omega)$. The magnetic field associated to this potential is a Dirac delta centered at $a$ and directed orthogonally to the plane. It determines the so called Aharonov--Bohm effect: a quantum particle moving in $\Omega\setminus\{a\}$ will be affected by the magnetic potential, although it remains in a region where the magnetic field is zero. This phenomenon can be simulated experimentally by the presence of a thin solenoid placed at $a$ and orthogonal to the plane (see \cite{AB}). We are concerned with the analysis of $K_\psi$--real solutions of the equation $H_{\mpAa,V}=0$, according to the following definition
\begin{defi}
We say that $U:\Omega \to \C$ is $K_\psi$--real if there exists a function $\psi \in C^1(\Omega\setminus\{a\},\C)$ with
\begin{equation}\label{hypotesis_psi}
|\psi|=1, \qquad \text{deg}(\psi, a)=2n+1, \ n\in\Z
\end{equation}
such that
\[
U = \psi \bar{U}.
\]
\end{defi}
This is a generalization of the notion of $K$--real function given in \cite{HHO,HHO1,BHH}, where the authors consider $\psi=e^{i\theta}$ ($\theta$ is the angular coordinate centered at the singular point). As the authors notice therein, this condition arises quite naturally in this context. In fact it is proved in \cite{HHO1} that the spectrum of $H_{\mpAa,V}$ (with homogeneous boundary conditions) consists of eigenvalues corresponding to $K_\psi$--real eigenfunctions. We are interested in non--homogeneous boundary conditions, in particular we will study the following systems
\begin{equation} \label{cauchy_pb_A}
\left\{ \begin{array}{ll}
(i\nabla +\mpAa)^2 U + V U = 0 \quad & \mathrm{in} \ \Omega\setminus \{a\}  \\
U= \Gamma & \mathrm{on} \ \partial \Omega,
\end{array} \right.
\end{equation}
where $a$ (the concentration point of the magnetic field) is intended as a parameter, whereas the boundary data $\Gamma\in W^{1,\infty}(\partial\Omega,\C)$ is fixed. A natural condition on the boundary trace is therefore
\begin{defi}
We say that $\Gamma:\partial \Omega \to \C$ is a $K_\psi$--real trace if there exists a function $\psi \in C^1(\partial\Omega,\C)$ satisfying (\ref{hypotesis_psi}), such that $\Gamma = \psi \bar{\Gamma}$.
\end{defi}
This assumption on the boundary trace implies that, for a suitable choice of the gauge in (\ref{potential}), the solution of (\ref{cauchy_pb_A}) is $K_\psi$--real. Therefore such a solution can be identified with a real function on the double covering of $\Omega\setminus \{a\}$. More precisely, it corresponds to an antisymmetric solution of a real elliptic equation on the double covering. Hence the nodal set
$$
\mathcal{N}(U) = \overline{\{x \ \in \ \Omega \ :\ U(x)=0\}}
$$
inherits some regularity properties, that is
\begin{teo} \label{theorem_nodal_set_U_1}
Let $\Gamma\in W^{1,\infty}(\partial\Omega,\C)$ be $K_\psi$--real and $a\in\Omega$ be fixed. For a suitable choice of the gauge in (\ref{potential}), the solution $U$ of \eqref{cauchy_pb_A} satisfies the following properties.
\begin{itemize}
\item [(i)] The nodal set of $U$ is nontrivial and consists of the union of regular arcs, having endpoints either at $\partial \Omega$, or at an interior singular point of $U$, or at $a$. Moreover there is at least one nodal line ending at $a$. 
\item [(ii)] There is an odd number of nodal lines ending at $a$. There is an even number of arcs meeting at interior singular points different from $a$.
\item [(iii)] Let $(r_a,\thetaa)$ be the polar coordinates centered at $a$. Then $U$ satisfies, for some odd $k \geq 1$, the asymptotic formula
\begin{equation} \label{sviluppo_U}
U(r_a,\theta_a)=e^{i \Theta(r_a,\theta_a)} \frac{{r_a}^\frac{k}{2}}{k}\left[ c_{k} \cos\left(\frac{k}{2}\theta_a\right) + d_{k} \sin\left(\frac{k}{2}\theta_a\right) \right] + o({r_a}^\frac{k}{2}),
\end{equation}
where $c_k^2+d_k^2\neq0$ and $e^{i \Theta}$ is a suitable complex phase. In particular, there is an odd number of nodal lines ending at $a$.
\end{itemize}
\end{teo}

Since we are mainly interested in triple partitions of the domain, which is the easiest case of multiple clustering, we shall require in addition that the boundary trace vanishes exactly three times on $\partial\Omega$. More precisely, we introduce the class
$$
\mathcal{G}=\{ \Gamma:\partial\Omega \to \C : \quad \Gamma \text{ is }K_\psi\text{--real} \quad \text{and} \quad
|\Gamma|=\sum_{i=1}^3\gamma_i \text{ with } (\gamma_i) \in g \},
$$
where
\begin{equation*}
g = \left\{ (\gamma_1,\gamma_2,\gamma_3) \text{ with } \gamma_i:\partial\Omega \to \R: \begin{array}{l}
				\gamma_i\cdot\gamma_j=0 \quad \textrm{for} \ i \neq j,\\
				\gamma_i \geq 0, \quad \gamma_i \in C^1(\overline{\{\gamma_i >0\}}), \quad i=1,2,3,\\
				\sum_{i=1}^3\gamma_i \ \textnormal{vanishes exactly three times on}\ \partial \Omega\\
\end{array}
\right\}.
\end{equation*}
In order to ensure that the nodal lines of $U$ do not dissect the domain in more than three parts, we will also impose the following coercivity condition on the scalar potential
\begin{equation} \label{hypothesis_V}
\int_\Omega \left( |\nabla\phi|^2+V\phi^2 \right) dx_1 dx_2 >0 \quad \forall \ \phi \in H^1_0(\Omega),
\end{equation}
then the following holds.
\begin{teo} \label{theorem_nodal_set_U_2}
Under the same assumptions of Theorem \ref{theorem_nodal_set_U_1}, suppose moreover that $\Gamma\in \mathcal{G}$ and that (\ref{hypothesis_V}) holds. Then the nodal set of $U$ consists of at most three arcs. If $\mathcal{N}(U)$ consists of one or two arcs, then the expansion \eqref{sviluppo_U} holds with $k=1$. Otherwise the three arcs intersect in exactly one point, that is $a$, and the expansion \eqref{sviluppo_U} holds with $k=3$.
\end{teo}
In this last case, the nodal set determines a partition of $\Omega$ into three parts: we will say that $a$ is a \textit{triple point} for $\Gamma$. Our aim is to understand the circumstances related to the occurrence of triple points, and in particular to investigate the variational characterizations of the triple point configuration. The first result involves the quadratic form associated to the operator
\begin{equation*}
Q_{\mpAa,V}(U)=\int_{\Omega}\left(\vert (i\nabla+\mpAa)U \vert^2+V\vert U\vert^2\right) dx_1 dx_2,
\end{equation*}
acting complex functions sharing the same boundary trace $\Gamma$. In this direction we prove
\begin{teo}[(Criticality)]\label{theorem_criticality}
Let $\Gamma\in\mathcal{G}$ and correspondingly let $\mpAa$ be a suitable choice of the gauge in (\ref{potential}). If the function
\[
a\mapsto \varphi(a)=\min\{Q_{\mpAa,V}(U)\;:\;U\in H^1(\Omega), \; U= \Gamma \;\mathrm{on} \; \partial \Omega\}.
\]
is differentiable\footnote{The differentiability of $\varphi(a)$ requires regularity of the boundary trace.} and (\ref{hypothesis_V}) holds, then the only critical points of $\varphi(a)$ are the triple points.
\end{teo}

Our next result is related to the conjecture proposed in \cite{BHH}, Section 6. The authors show that the second eigenfunction of $H_{\mpAa,0}$ with homegenous Dirichlet boundary conditions, whenever it exhibits triple point configuration, turns to be a reasonable candidate for the spectral minimal 3--partition. This is supported by numerical simulations (see also \cite{BHV} for related problems). Similarly, we succeed in characterizing the triple point configuration as minimal 3--partition, with respect to the real energy functional associated to the equation. Given $\Gamma \in \mathcal{G}$, let $(\gamma_i)\in g$ be the associated real trace. For $i=1,2,3$, we define a set function $J_i(\omega_i)$, acting on open sets $\omega_i\subset \Omega$, satisfying $\textnormal{supp}(\gamma_i)\subset \partial \Omega\cap\partial \omega_i$, in the following way
\begin{equation*}
J_i(\omega_i)=\inf\left\{ \int_\Omega \left( |\nabla u|^2 + V u^2 \right) dx_1 dx_2 :\begin{array}{l}
                                    u \in H^1(\Omega), \ u=\gamma_i \ \textrm{on} \  \partial \Omega\cap\partial\omega_i \\
				    u\geq 0, \ u=0 \ \textrm{in } \Omega \setminus\omega_i
                                    \end{array}
\right\}.
\end{equation*}
Starting from the results in \cite{ctv3}, we can prove
\begin{teo}[(Optimal partition)]\label{theorem_minimal_partition}
Assume that \eqref{hypothesis_V} holds. Suppose moreover that $\Gamma \in \mathcal{G}$ admits triple point $a$, so that $\Omega \setminus \mathcal{N}(U)$ has three connected components. Then the connected components are solution of the optimal partition problem
\begin{equation}\label{minimum_pb}
\inf\left\{ \sum_{i=1}^3 J_i(\omega_i) :\begin{array}{l}
                                    \omega_i \textnormal{ open},\ \textnormal{supp}(\gamma_i)\subset \partial \omega_i \\
				    \cup_{i=1}^3\overline{\omega_i}=\overline{\Omega}, \ \omega_i \cap \omega_j=\emptyset, 
\ i\neq j
                                    \end{array}
\right\}.
\end{equation}
\end{teo}

This variational characterization is achieved through a uniqueness result. For a given boundary trace, we analyze the set of systems (\ref{cauchy_pb_A}) as $a$ varies in $\Omega$ and prove uniqueness of the triple point configuration with respect to the parameter $a$.
\begin{teo}[(Global uniqueness)] \label{theorem_global_uniqueness}
If (\ref{hypothesis_V}) holds, then every $\Gamma \in \mathcal{G}$ admits at most one triple point. Moreover the set of boundary traces which admit a triple point is open and dense in $\mathcal{G}$ (with respect to the $L^\infty$--norm). Finally, the position of the triple point depends continuously on the $L^\infty$--norm of the boundary data.
\end{teo}

We conclude the analysis by showing that the $C^{1}$--norm of the nodal lines depends continuously on the $L^\infty$--norm of the boundary data.

Our proof of Theorems \ref{theorem_minimal_partition}, \ref{theorem_global_uniqueness} is strongly related to some previous works by Conti, Terracini and Verzini. On the other hand, we also generalize some results therein, providing
\begin{teo}\label{theorem_uniqueness_Sgamma}
Assume that (\ref{hypothesis_V}) holds. For every $(\gamma_i) \in g$, define the functional class
\begin{equation*}
\Sgamma=\left\{ (u_1,u_2,u_3) \in (H^1(\Omega))^3: \begin{array}{l}
                                    \ui \geq0, \ \ui=\gamma_i \ \textrm{on} \  \partial \Omega \\
				    \ui\cdot\uj=0 \ \textrm{a.e.} \ x\in \Omega, \ \textrm{for} \ i \neq j \\
				    -\Delta\ui +V\ui \leq 0, \quad  -\Delta \hatui +V \hatui \geq0 \quad i=1,2,3 
\end{array}
\right\}
\end{equation*}
where $\hatui := \ui-\sum_{j \neq i} \uj$. Then
\begin{itemize}
 \item[(i)] $\Sgamma$ consists of exactly one element $(u_1,u_2,u_3)$;
\item[(ii)] the open sets $\omega_i=\{u_i>0\}$ are the unique solution of (\ref{minimum_pb}).
\end{itemize}
\end{teo}
We wish to mention that this functional class also contains the limiting configurations of the solutions of a competition--diffusion system as the competition parameter tends to infinity (see \cite{ctv4}). Theorem \ref{theorem_global_uniqueness} also establishes a variational characterization of these limiting functions, which is quite surprising since the initial system is not variational (see Section \ref{section_energy_minimizing_partitions} for further details).


The paper is organized as follows. In Section \ref{sec_gauge} we introduce the notion of classical gauge invariance and show that equation \eqref{cauchy_pb_A} is equivalent to a real elliptic equation on the double covering of $\Omega\setminus\{a\}$. For this part we mainly refer to \cite{HHO1}, where the homogeneous Cauchy problem is studied. In Section \ref{section_nodal_set} we perform the analysis of the nodal lines and prove Theorems \ref{theorem_nodal_set_U_1}, \ref{theorem_nodal_set_U_2}. This is based on the well known properties of the nodal set of solutions to real elliptic equations in the plane; we address the reader to the classical works by Alessandrini \cite{AL} and Hartman and Wintner \cite{HW}. As it concerns the non planar case, we refer to \cite{HL,HHHN} and references therein. We also wish to mention some results concerning Schr\"odinger operators with singular potentials, such as \cite{HH,HHN3,HHN4}. For the specific case of Aharonov--Bohm type magnetic potentials, few is known because of the very strong singularity; the recent paper \cite{FFT} provides regularity results for a large class of equations including \eqref{cauchy_pb_A}. Section \ref{section_local_uniqueness} contains the technical part of the work. First of all we prove Theorem \ref{theorem_criticality} and then extend the result to the Rayleigh quotient, providing Theorem \ref{theorem_eigenvalues}. Finally we prove a local uniqueness result, that is a local version of Theorem \ref{theorem_global_uniqueness}. All the remaining stated theorems are proved in Section \ref{section_energy_minimizing_partitions} which, as we mentioned, is strongly based on the previous works by Conti, Terracini and Verzini \cite{ctv2,ctv4,ctv3,ctv5}. Finally, in Section \ref{section_continuous_dependence}, we prove the continuous dependence of the nodal lines with respect to the data.

\noindent{\bf Acknowledgments} The first author would like to thank Prof. M. Ramos for his kind invitation and financial support during the final typing of the work.

\section{Preliminaries} \label{section_schrodinger_aharonov_bohm}
In this section we fix the main assumptions and notations and we introduce the notion of solution of the equation $H_{\mpAa,V}=0$.

As it concerns the domain $\Omega$, in the following we will always work in the open unit disk of $\R^2$, which will be denoted by $D$. This is not restrictive thanks to the Riemann mapping theorem (see for example \cite{GK}, Theorem 6.42) and due to the conformal equivariance of the problem, which will be proved in Section \ref{sec_gauge}. Given a point $a \in D$ we denote by $\Da$ the open set $D\setminus \{a\}$.

The regularity assumptions on the physical quantities, which will hold throughout the paper, are the following
\[
\mpAa \in L^1(D)\cap C^1(\overline{D}\setminus\{a\}) \qquad\text{and}\qquad V \in W^{1,\infty}(D).
\]
Notice that the magnetic potential associated to the Aharonov--Bohm effect presents a strong singularity at $a$, in particular it is not in $L^2(D)$, hence this assumption on the magnetic potential is quite natural in this context. The magnetic field associated to (\ref{potential}) satisfies
\begin{equation} \label{magnetic_field}
\mfB_a= \curl \mpAa = (2n+1) \pi \delta_a \textbf{k} \quad \textrm{in} \ \Omega, \quad n\in\Z,
\end{equation}
where $\delta_a$ is the Dirac delta centered at $a$ and $\textbf{k}$ is the unitary vector orthogonal to the plane. It will be usefull to notice that this can equivalently substituted with the condition
\begin{equation}\label{HA1}
\curl \ \mpAa =0 \quad \mathrm{in} \ \Da,
\end{equation}
together with an additional assumption on the normalized circulation
\begin{equation}\label{HA2}
\frac{1}{2 \pi} \oint_\sigma \mpAa \cdot \bfdx =\frac{2n+1}{2}, \quad n \in \Z
\end{equation}
for every closed path $\sigma$ which winds once around the pole.

Let us now introduce the space $\HAD$ of the solutions of (\ref{cauchy_pb_A}). We refer to \cite{MR} for a complete review on magnetic Schr\"odinger operators and to \cite{MOR} for the specific case of the A--B effect. We recall that the operator acts on complex valued functions as
\begin{equation*}
H_{\mpAa,V}U =- \Delta U+i \div (\mpAa U) + i \mpAa \cdot \nabla U+|\mpAa|^2 U + V U.
\end{equation*}
Following \cite{FFT}, we define $\HAD$ as the completion of
\[
\mathcal{H}^0(D)=\{ U\in C^\infty(D,\C): \text{ $U$ vanishes in a neighborhood of } a \}
\]
with respect to the norm
\[
\| U\|_{\HAD} = \left( \| (i\nabla+\mpAa)U\|^2_{L^2(D,\C^2)}+\|U\|^2_{L^2(D,\C)} \right)^{1/2}.
\]
It is proved in \cite{LaWe} that a magnetic Hardy inequality holds in dimension two, whenever the circulation of $\mpAa$ is not integer. As pointed out in \cite{MOR}, such inequality holds also for functions that are not compactly supported. More precisely, for every $U\in\HAD$ it holds
\[
\int_D \frac{|U|^2}{|x-a|^2}dx_1 dx_2 \leq 4 \int_D |(i\nabla+\mpAa)U|^2 dx_1dx_2.
\]
This implies that $\HAD$ is a subset of the usual Sobolev space $H^1(D,\C)$ and we can give the following equivalent characterization.

\begin{lemma}
Let $\mpAa$ satisfy \eqref{potential}, then
\[
\HAD =\left\{ U\in H^1(D,\C) : \ \frac{U}{|x-a|}\in L^2(D,\C) \right\}.
\]
\end{lemma}
\begin{proof}
Assume first that $U\in H^1(D,\C)$ and $U/|x-a|\in L^2(D,\C)$, then using expression (\ref{potential}) for the potential we obtain
\begin{eqnarray*}
\|U\|^2_{\HAD} &\leq& 2 \|U\|^2_{H^1(D,\C)} +2 \int_D|\mpAa U|^2 dx_1 dx_2 \\
&\leq& 2 \|U\|^2_{H^1(D,\C)} + 2 \|(x-a)\mpAa\|_{L^\infty(D)}^2 \int_D \frac{|U|^2}{|x-a|^2} dx_1dx_2 < \infty.
\end{eqnarray*}
On the other hand if $U\in \HAD$ then
\begin{eqnarray*}
\|U\|^2_{H^1(D,\C)}&\leq& \| (i\nabla+\mpAa)U\|^2_{L^2(D,\C^2)} +\| \mpAa U \|^2_{L^2(D,\C^2)}+\| U\|^2_{L^2(D,\C)} \\
&\leq& \| U \|_{\HAD}^2 +\| (x-a)\mpAa \|_{L^\infty(D)}^2 \int_D \frac{|U|^2}{|x-a|^2} dx_1dx_2 \\
&\leq& C \| U \|_{\HAD}^2,
\end{eqnarray*}
where we used the magnetic Hardy inequality in the last step.
\end{proof}
The previous lemma implies in particular that, whenever (\ref{hypothesis_V}) holds, the quadratic form $Q_{\mpAa,V}$ associated to the operator is definite positive.

Given a boundary data $\Gamma \in H^{1/2}(\partial D,\C)$, we say that $U\in \HAD$ is a solution of \eqref{cauchy_pb_A} if the following integral equality holds for every $\phi \in \HAD$
$$
\int_D [ (i\nabla +\mpAa) U \overline{(i\nabla +\mpAa)\phi} +V U \overline{\phi} ] dx_1 dx_2 
+ i \int_{\partial D} [\Gamma \overline{(i\nabla +\mpAa)\phi}\cdot \nu + (i\nabla +\mpAa)U\cdot \nu \overline{\phi}] d\sigma =0.
$$
Notice that $\mathcal{G}\subset W^{1,\infty}(\partial\Omega,\C)$ and in fact we will need to work with this regularity on the boundary trace.

\section{Gauge invariance for A--B potentials with semi--integer circulation} \label{sec_gauge}

In this section we shall present in detail a result contained in \cite{HHO1}, related to the gauge invariance property of magnetic operators of A--B type, having semi-integer circulation. Our aim is to generalize it to the non-homogeneous Cauchy problem, in the case of $K_\psi$--real boundary traces; this will be done in Proposition \ref{prop_U_u} by means of a suitable choice of the gauge in (\ref{potential}). Throughout this section $a\in D$ is fixed. Let us start introducing the general notion of gauge equivalence for magnetic Schr\"odinger operators.

Let $\Omega \subset \R^2$ be bounded domain, $\Omegatilde$ be a covering manifold and $\Pi: \Omegatilde \to \Omega$ be the associated projection map. We endow $\Omegatilde$ with the locally flat metric obtained by lifting the Euclidean metric of $\Omega$, in such a way that $\Pi$ is a local isometry. Since the differential and integral operators on $\Omegatilde$ coincide locally with the usual ones, we will denote them with the same symbol.

\begin{defi} \label{lift} For a function $f:\Omega \to \C$ we define the lifted function $\tilde{f}:\Omegatilde \to \C$ as
$$
\tilde{f}=f \circ \Pi.
$$
For a path $\sigma: [0,1] \to \Omega$ and a point $p \in \Omegatilde$ such that $\Pi(p)=\sigma(0)$ let $\sigmatilde : [0,1]\to \Omegatilde$ denote the unique lifted path such that $\sigmatilde(0)=p$ and 
$$
\sigma=\Pi \circ \sigmatilde.
$$
\end{defi}

\begin{lemma} \label{lemma_gauge_equivalence}
Let $\mpAa,\mpAa'$ satisfy \eqref{HA1} in $\Da$ and assume moreover that
$$
\frac{1}{2 \pi} \oint_\sigma (\mpAa'-\mpAa) \cdot \bfdx \ \in \ \Z,
$$
for every closed path $\sigma$ in $\Da$. Then, denoting $\tilde{D}_a$ the universal covering manifold of $\Da$, there exists $\Theta \in C^2(\tilde{D}_a)$ such that\footnote{Here and in the following we may write, with some abuse of notation, $\nabla\Theta$ instead of $\Pi(\nabla\Theta)$.} $\nabla\Theta\in C^1(\Da)$ and
\begin{equation}\label{eq:gauge_equivalence}
\mpAa'=\mpAa+\nabla \Theta \text{ in } \Da, \qquad
e^{i\Theta}\in C^2(\Da).
\end{equation}
\end{lemma}

\begin{proof} Due to the fact that $\tilde{D}_a$ is simply connected and $\mpAa,\mpAa'$ satisfy \eqref{HA1}, it holds
$
\oint_{\sigmatilde} (\mpAtilde'-\mpAtilde) \cdot \bfdx=0,
$
for every closed path $\sigmatilde$ in $\tilde{D}_a$. Therefore there exists $\Theta \in C^2(\tilde{D}_a)$ such that $\mpAtilde_a'=\mpAtilde_a+\nabla \Theta$ in  $\tilde{D}_a$. Consider now any two points $p,p' \in \tilde{D}_a$ such that $\Pi(p)=\Pi(p')$. Then for every path $\sigmatilde$ on $\tilde{D}_a$ connecting $p$ to $p'$ we have
\begin{eqnarray}
\Theta(p)-\Theta(p')&=&\int_{\sigmatilde} \nabla \Theta \cdot \bfdx 
= \int_{\sigmatilde} (\mpAtilde_a'-\mpAtilde_a) \cdot \bfdx \nonumber \\
&=& \int_{\sigma} (\mpAa'-\mpAa) \cdot \bfdx = 2\pi n, \nonumber
\end{eqnarray}
for some $n \in \Z$. Therefore both $e^{i\Theta}$ and $\nabla\Theta$ are well defined in $\Da$.
\end{proof}

Motivated by this lemma, we give the following definition of gauge equivalence.

\begin{defi}
We say that $\mpAa,\mpAa'$ satisfying \eqref{HA1} are gauge equivalent in $\Da$ if there exists $\Theta \in C^2(\tilde{D}_a)$ such that
\eqref{eq:gauge_equivalence} holds (we can also say that the operators $H_{\mpAa,V},H_{\mpAa',V}$ are gauge equivalent).
\end{defi}

If in particular the circulation of $\mpAa-\mpAa'$ vanishes, then $\Theta\in C^2(D)$; this implies that every magnetic potential satisfying (\ref{magnetic_field}) has the form (\ref{potential}). The next result shows that the solutions of gauge equivalent operators only differ by a complex phase, hence in particular they share the same nodal set.

\begin{lemma} \label{lemma_gauge_invariance}
Let $\mpAa,\mpAa'$ as in the previous lemma. If $U\in\HAD$ is a solution of $H_{\mpAa,V} U=0$ in $\Da$, with Dirichlet boundary conditions, then $U'=e^{i \Theta} U$ is solution of $H_{\mpAa',V} U'=0$ in $\Da$.
\end{lemma}

\begin{proof}
Let $\phi \in C_0^{\infty}(\Da)$ be a test function and $\Theta \in C^2(\tilde{\Da})$ as in the previous lemma, then a direct calculation shows that $e^{i \Theta}(i\nabla+\mpAa)^2\phi=(i\nabla+\mpAa')^2(e^{i \Theta}\phi)$. By multiplying the equation by $e^{i \Theta}e^{-i \Theta}$, we obtain
$$
0= \int_D U [\overline{(i\nabla +\mpAa)^2\phi +V\phi} ] dx_1 dx_2 =
\int_D U' [\overline{(i\nabla +\mpAa')^2(e^{i \Theta}\phi) +Ve^{i \Theta}\phi} ] dx_1 dx_2.
$$
Since $e^{i\Theta} \in C^2(\Da)$, the result follows by density.
\end{proof}

\subsection{The twofold covering manifold}
As a particular case of Lemma \ref{lemma_gauge_equivalence} we infer that, whenever the circulation of $\mpAa$ is an integer, the magnetic operator $H_{\mpAa,V}$ is gauge equivalent to the elliptic operator $H_{0,V}=-\Delta+V$. Let us now turn to consider potentials satisfying (\ref{HA2}). In this case we can still relate $H_{\mpAa,V}$ and $H_{0,V}$, provided we replace the domain $\Da$ with its twofold covering manifold. This result is proved in \cite{HHO} for the Dirichlet and Neumann homogeneous cases. We point out that the operators are \textit{not} unitarily equivalent since, as we are going to see, there is a one-to-one correspondence between the eigenfunctions of $H_{\mpAa,V}$ and the antisymmetric eigenfunctions of $H_{0,V}$.

\begin{defi}
The twofold covering manifold of $\Da$ is the following subset of $\C^2$ 
\begin{equation*}
\Sigmaa= \{ (x,y) \ \in \ \C^2 : y^2=x-a, \ x \in \Da  \},
\end{equation*}
endowed with the locally flat Euclidean metric. We will denote
$$
\Pi_x :(x,y) \mapsto x \qquad \Pi_y :(x,y) \mapsto y,
$$
the two projections naturally defined on $\Sigmaa$.
\end{defi}

The following proposition describes more explicitly the locally flat metric considered on $\Sigmaa$.

\begin{prop} \label{prop_sigma}
There exists a global chart of $\Sigmaa$ which coincides locally with $\Pi_x$. In particular it induces on $\Sigmaa$ a locally flat Euclidean metric.
\end{prop}

\begin{proof} Let $(r_a,\thetaa)$ be the polar coordinates centered at $a$, in such a way that
$$
y^2(x)=x-a=r_a(x) e^{i \theta_a(x)}.
$$
Then we can define the following parametrization: 
\begin{eqnarray}
\Phi : [0,4 \pi) \times [0,1) &\to& \Sigmaa \cup \{(a,0)\} \nonumber \\
(\theta,r) &\mapsto& \left( re^{i\theta}, \sqrt{r_a(r,\theta)} e^{i \frac{\theta_a (r, \theta)}{2}} \right). \nonumber
\end{eqnarray} 
The function $\Phi$ is bijective on $[0,4 \pi) \times [0,1) \setminus \Phi^{-1}(a,0)$, therefore its inverse $\Phi^{-1}$ is well defined on this domain, and it is the desired chart.
\end{proof}

\begin{defi}\label{def_polar_coordinates}
We shall use the following notation for the polar coordinates on $\Sigmaa$
\[
x=r e^{i \theta}\qquad\text{and}\qquad x-a=r_a e^{i \theta_a},
\]
while
\[
y=\rho e^{i \phi} \qquad\text{with the relation}\qquad \rho=\sqrt{r_a}, \ \phi=\frac{\theta_a}{2}.
\]
\end{defi}

\begin{rem} \label{discontinuity_theta_a}
In the definition of the angle $\theta_a$ we usually consider it a discontinuous function on a horizontal segment starting at the point $a$. Nevertheless we can decide to move the discontinuity without altering the previous construction. In the future analysis in particular it will be useful to consider $\theta_a$ discontinuous on two adjacent segments (when $a\neq 0$): the segment connecting the origin with $a$ and the segment connecting the origin with a point $x_0 \in \partial \Omega$.
\end{rem}

\begin{defi} \label{symmetry_map_G}
On the twofold covering manifold we define a symmetry map $G: \Sigmaa \to \Sigmaa$, which associates to every $(x,y)$ the unique $G(x,y)$ such that $\Pi_x((x,y))=\Pi_x(G(x,y))$, that is $G(x,y):=(x, -y)$. We say that a function $f:\Sigmaa \to \C$ is symmetric if $f(G(x,y))=f(x,y), \ \forall (x,y) \in \Sigmaa$, and antisymmetric if $f(G(x,y))=-f(x,y), \ \forall (x,y) \in \Sigmaa$.
\end{defi}

Every function $f$ defined in $\Da$ can be lifted on $\Sigmaa$ as described in Definition \ref{lift}, by means of the projection $\Pi_x$. Notice that $\tilde{f}$ is always symmetric in the sense of the preceding definition. 

\begin{lemma}\label{lemma_psi_square}
Let $\Gamma$ be $K_\psi$--real, then there exists an antisymmetric function, that we denote by $\psi^{1/2}$, such that
\[
\psi^{1/2} \in C^1(\partial \Sigmaa,\C) \quad\text{and}\quad
\Pi_x\left[(\psi^{1/2})^2\right]=\psi,
\]
where $(\cdot)^2$ denotes the complex square. In particular $\psi^{-1/2}\Gammatilde$ is real valued and antisymmetric.
\end{lemma}

\begin{proof}
Let $\tilde{\psi}$ be the symmetric lifting of $\psi$ on $\partial\Sigmaa$, then
\[
|\tilde{\psi}|=1 \quad\text{ and }\quad \text{deg}\ \tilde{\psi}=2(2n+1).
\]
Therefore the composition with the complex square root function $\psi^{1/2}:= \sqrt{\tilde{\psi}}$ is well defined on $\partial \Sigmaa$, providing the first part of the statement. The second part comes from the fact that, by definition, $\Gamma=\psi\bar{\Gamma}$ and finally the complex square root is antisymmetric on $\partial \Sigmaa$.
\end{proof}

The next result shows that $H_{\mpAa,V}$ is gauge equivalent to $H_{0,V}$ in $\Sigmaa$. Given a boundary trace, we select the gauge in such a way to obtain $K_\psi$--real solutions for every choice of the position of the singularity $a$.

\begin{lemma}\label{lemma_potenziale2}
Let $\Gamma$ be $K_\psi$--real and let $\psi^{1/2}$ as in the previous lemma. Then there exists $\Thetaa\in C^2(\tilde{\Da})$ such that\footnote{Here $\tilde{\Da}$ is the universal covering of $\Da$.}
\begin{itemize}
\item[(i)] $e^{i\Thetaa}\in C^2(\Sigmaa)$ is antisymmetric and $e^{i\Thetaa}|_{\partial\Sigmaa}=\psi^{1/2}$;
\item[(ii)] the potential $\mpAa:=\nabla\Thetaa$ satisfies (\ref{magnetic_field}), hence in particular $H_{\mpAa,V}$ is gauge equivalent to $H_{0,V}$ in $\Sigmaa$.
\end{itemize}
\end{lemma}

\begin{proof}
By the definition of $\psi$, it is possible to choose the discontinuity of $\thetaa$ as in Remark \ref{discontinuity_theta_a} in such a way that $-\frac{i}{2}\log\psi-\frac{2n+1}{2}\thetaa$ is continuous on $\partial D$. Hence we can consider its harmonic extension on the disk
\begin{equation} \label{sys:section_gauge}
\left\{ \begin{array}{ll}
-\Delta \Phi = 0 \quad & \mathrm{in} \ D  \\
\Phi= -\frac{i}{2}\log\psi-\frac{2n+1}{2}\thetaa & \mathrm{on} \ \partial D.
\end{array} \right.
\end{equation}
Then the desired potential is the function $\Thetaa: \tilde{\Da} \to \R$ defined by
\[
\Thetaa:=\frac{2n+1}{2}\thetaa+\tilde{\Phi}.
\]
Clearly $e^{i\Thetaa}$ is well defined on $\Sigmaa$ and $e^{i\Thetaa}|_{\partial\Sigmaa}=\psi^{1/2}$. Moreover $\nabla \Thetaa$ can be projected in $\Da$ and $\mpAa:=\nabla\Thetaa$ has the form (\ref{potential}), with $\Phi$ defined in (\ref{sys:section_gauge}). Let us show that $e^{i\Thetaa}$ is antisymmetric. Fix $(x,y) \in \Sigmaa$ and let $\sigmatilde :[0,1] \to \Sigmaa$ be a path which joins $(x,y)$ to $G(x,y)$, then using the notations of Definition \ref{lift}, it holds
$$
\frac{1}{2 \pi} \oint_{\sigmatilde} \mpAtilde \cdot \bfdx =
\frac{1}{2 \pi} \oint_{\sigma} \mpAa \cdot \bfdx =\frac{2n+1}{2}, \quad n \in \Z.
$$
Therefore:
$$
\Thetaa(G(x,y)) -\Thetaa(x,y) = \int_{\sigmatilde} \nabla \Thetaa \cdot \bfdx = (2n+1)\pi,
$$
and hence $e^{i\Thetaa(G(x,y))}=-e^{i\Thetaa(x,y)}$. Finally let $\sigmatilde$ be a closed path in $\Sigmaa$, then
$$
\frac{1}{2 \pi} \oint_{\sigmatilde} \mpAtilde_a \cdot \bfdx
= \frac{1}{2 \pi} \oint_{\sigma} \mpAa \cdot \bfdx  = m \ \in \ \Z,
$$
since $\sigma$ always turns an even number of times around the singularity. Hence the gauge equivalence comes from Lemma \ref{lemma_gauge_equivalence}.
\end{proof}

We can finally prove the existence of a bijection between the solutions of (\ref{cauchy_pb_A}) and the antisymmetric solutions of a real elliptic problem on the twofold covering manifold.

\begin{prop} \label{prop_U_u}
Let $\Gamma$, $\Thetaa$ and $\mpAa$ be as in the previous lemma, and denote by $U$ the corresponding solution of (\ref{cauchy_pb_A}). Then the function
\begin{equation} \label{U_u}
u(x,y):= e^{-i \Thetaa(x,y)} \tilde{U}(x,y).
\end{equation}
is antisymmetric, real valued and solves
\begin{equation} \label{cauchy_pb_0}
\left\{ \begin{array}{ll}
-\Delta u + \tilde{V} u = 0 \quad & \mathrm{in} \ \Sigmaa  \\
u= \gamma & \mathrm{on} \ \partial \Sigmaa,
\end{array} \right.
\end{equation}
where $\gamma=\psi^{-1/2}\Gammatilde$.
\end{prop}

\begin{proof}
We proceed as in Lemma \ref{lemma_gauge_invariance}, but taking into account the boundary trace. If $\phi\in \HAD$ is any test function it holds
\begin{eqnarray*}
0 &=& \int_{\Sigmaa} \tilde U [\overline{(i\nabla +\tilde \mpAa)^2\phi +V\phi} ] dx_1 dx_2 
+ i \int_{\partial \Sigmaa} [\tilde \Gamma \overline{(i\nabla +\tilde \mpAa)\phi}\cdot \nu 
+ (i\nabla +\tilde \mpAa)\tilde U\cdot \nu \overline{\phi}] d\sigma \\
&=&\int_{\Sigmaa} u [\overline{(-\Delta)(e^{-i \Thetaa}\phi) +\tilde V e^{-i \Thetaa}\phi} ] dx_1 dx_2
+ i \int_{\partial \Sigmaa} [\gamma \overline{(i\nabla) (e^{-i \Thetaa}\phi)}\cdot \nu
+ i\nabla u\cdot \nu \overline{(e^{-i \Thetaa}\phi)}] d\sigma. \\
\end{eqnarray*}
Hence for every real valued test function $\psi$ we have
\[
\int_{\Sigmaa}(-\Delta \psi+\tilde V \psi) u dx_1 dx_2 
+\int_{\partial \Sigmaa}(\gamma \nabla\psi\cdot \nu - \psi\nabla u \cdot \nu)d\sigma=0,
\]
which is the weak form of \eqref{cauchy_pb_0}. By Lemma \ref{lemma_psi_square}, $\gamma$ is real valued, hence so is $u$. Moreover $u$ is the product of an antisymmetric function times a symmetric one, hence it is antisymmetric.
\end{proof}

From now on we will always make this choice of the gauge, in such a way to consider only $K_\psi$--real solutions.

\subsection{Related real elliptic problems} \label{subsec_related_real_elliptic_problems}
In this subsection we will prove the existence of a bijection between the solutions of (\ref{cauchy_pb_A}) and the antisymmetric solutions of a real elliptic equation in a bounded subset of $\R^2$. This is performed in two different ways. In the following lemma we simply apply the projection $\Pi_y$ to the function $u$ defined in (\ref{U_u}). We obtain a real valued function which will be suitable for the local analysis that we perform in Section \ref{section_nodal_set}.

Here and in the following we will often make the identification $\R^2 \backsimeq \C$, writing $x=(x_1,x_2)=x_1+ix_2$. We shall use the following standard notation for the complex derivative
$$
\parder{}{x}=\frac{1}{2}\left(\parder{}{x_1}-i \parder{}{x_2}\right) \quad 
\parder{}{\bar{x}}=\frac{1}{2}\left(\parder{}{x_1}+i \parder{}{x_2}\right),
$$
where $\parder{}{x_i}$ denotes partial derivative.

\begin{lemma} \label{lemma_definition_uuno}
In the same assumptions and notations of Proposition \ref{prop_U_u}, the function $\uuno(y):=u \circ \Pi_y^{-1}(y)$ is an odd solution of the real elliptic equation
\begin{equation*}
\left\{ \begin{array}{ll}
-\Delta \uuno + \Vuno \uuno = 0 \quad & \mathrm{in} \ \Omegaa  \\
\uuno= \gammauno & \mathrm{on} \ \partial \Omegaa,
\end{array} \right.
\end{equation*}
where $\Vuno(y)=4|y|^2 \tilde{V} \circ \Pi_y^{-1}(y)$ and $\gammauno(y)=\gamma \circ \Pi_y^{-1}(y)$. Moreover $\uuno \in C^2_{loc}\cap W^{1,p}(\Omegaa)$, for every $p<+\infty$ and
$$
Q_{\mpAa,V}(U)=\frac{1}{2} \int_{\Omegaa} (|\nabla\uuno|^2+\Vuno(\uuno)^2)dy_1dy_2.
$$
\end{lemma}

\begin{proof}
Notice that the projection $\Pi_y$ is a global diffeomorphism of $\Sigmaa$ onto its image, hence its inverse is well defined
\begin{eqnarray}
\Pi_y^{-1}: & \Pi_y(\Sigmaa) & \to \Sigmaa \nonumber \\
& y & \mapsto (y^2+a, y). \nonumber
\end{eqnarray}
It is immediate to see that $\Pi_y^{-1}$ is conformal on its domain with respect to the metric defined on $\Sigmaa$ (by composing with the chart $\Phi$ defined in Proposition \ref{prop_sigma}); its conformal factor is
$$
\left| \parder{(\Pi_y^{-1})}{y}(y) \right|^2=4|y|^2.
$$
Now, being a bounded map, it admits a conformal extension in $\Pi_y(\Sigmaa)\cup\{(0,0)\}$, see for example \cite{GK}, Proposition 4.3.3. Hence the statement comes from the well known properties of conformal maps.
\end{proof}

\begin{defi}
In the following we will denote by $\Omegaa$ the set $\Pi_y(\Sigmaa)\cup\{(0,0)\}$. Notice that the singular point $a \in D$ corresponds to the origin in $\Omegaa$.
\end{defi}

In order to obtain a function defined in the unit disk (for every position of the singularity $a$) let us compose with a M\"obius transformation. This will be more appropriate for the analysis in Section \ref{section_local_uniqueness}, where the parameter $a$ varies.

\begin{lemma} \label{lemma_definition_udue}
Let $\uuno:\Omegaa \to \R$ as in the previous lemma. There exists a conformal map $T_a':D \to \Omegaa$ such that $\udue(y):=\uuno\circ T_a'(y)$ satisfies
\begin{equation*}
\left\{ \begin{array}{ll}
-\Delta \udue + \Vdue \udue = 0 \quad & \mathrm{in} \ D  \\
\udue= \gammadue & \mathrm{on} \ \partial D,
\end{array} \right.
\end{equation*}
where $\Vdue(y)=|\parder{T_a'}{y}(y)|^2 \Vuno \circ T_a'(y)$ and $\gammadue(y)=\gammauno \circ T_a'(y)$.
\end{lemma}

\begin{proof}
Proceeding as in \cite{ctv5} we consider the M\"obius transformation:
\begin{equation}\label{eq:moebius}
T_a:\overline{D}\longrightarrow\overline{D},\qquad
T_a(x)=\dfrac{x+a}{\bar a x+1}.
\end{equation}
It is well known that $T_a$ is a conformal map, such that $T_a(\partial D)=\partial D$ and $T_a({0})=a$. Let now $\tilde{T}_a(x,y):\Sigma_0\to\Sigmaa$ be the lifting of $T_a$. More precisely, if we denote for the moment $re^{i\theta}:=T_a(x)-a$, we have $\tilde{T}_a(x,y)=(re^{i\theta}+a,\sqrt{r}e^{i\frac{\theta}{2}})$. It only remains to prove that the map $T_a':D\to \Omega_a$, defined by $T_a'=\Pi_y\circ\tilde{T}_a\circ\Pi_y^{-1}$, is conformal. Indeed it is clearly conformal outside the origin, since the complex square root is conformal on the twofold covering manifold $\Sigmaa$, moreover, being bounded, it admits a conformal extension at the origin.
\end{proof}

\begin{rem}\label{remark_Vdue}
In the previous lemma we have equivalently
$$
\Vdue(y)=4|y|^2 \left|\parder{T_a}{x}(y^2)\right|^2 V \circ T_a(y^2), \qquad
\gammadue(y)=\gamma \circ T_a(y^2).
$$
\end{rem}

We end this section recalling a complex formulation of Green's theorem that we will need later.

\begin{lemma}\label{lemma_green}
Let $\Omega \subset \C$ be a regular domain and $\Phi,\Psi \in C^1(\overline \Omega,\C)$. Let $\mpA\in C^1(\overline \Omega)$ be a vector potential, such that $\nabla\Theta=\mpA$. Then it holds
\begin{multline*}
\int_{\Omega} \Psi (i \nabla +\mpA )^2 \Phi \ d x_1 d x_2 = 
2i\int_{\partial \Omega} \Psi \left(\parder{ }{x}-i \parder{\Theta}{x} \right) \Phi \ dx + \\
+ 4 \int_\Omega \left(\parder{ }{x}-i \parder{\Theta}{x} \right) \Phi \cdot \left(\parder{ }{\bar{x}}+i \parder{\Theta}{\bar{x}} \right) \Psi \ d x_1d x_2.
\end{multline*}
\end{lemma}

\begin{proof}
It is sufficient to apply the following complex formulation of Green's formula
\[
\int_{\Omega} \parder{F}{\xbar} dx_1 dx_2 =-\frac{i}{2} \int_{\partial \Omega} F dx, \quad \text{ with } \quad
F= -4 \parder{}{x}\left( e^{-i\Theta}\Phi \right) e^{i\Theta}\Psi,
\]
see for example \cite{GK}, Appendix A.
\end{proof}

\section{Properties of the nodal set}\label{section_nodal_set}
The aim of this section is the analysis of the nodal set of the solutions of (\ref{cauchy_pb_A}) and in particular the proof of Theorems \ref{theorem_nodal_set_U_1} and \ref{theorem_nodal_set_U_2}. Let us start recalling some known properties of the nodal set and singular points of solutions of real elliptic equations of the following kind
\begin{eqnarray}\label{eq:elliptic_generic}
\left\{ \begin{array}{ll}
-\Delta f + V f = 0 \quad & \mathrm{in} \ D  \\
f= \gamma & \mathrm{on} \ \partial D,
\end{array} \right.
\end{eqnarray}
with $V \in L^{\infty}(D)$, $\gamma \in W^{1,\infty}(\partial D)$. By standard regularity results and Sobolev imbedding, $f \in C^{1,\alpha}_{loc}(D) \cap W^{1,p}(D), \ \forall \alpha \in (0,1), p<+\infty$.

\begin{defi} We say that $y_0 \in \mathcal{N}(f)$ is a singular point if $\nabla f(y_0)=0$. We say that it is a zero of order (or multiplicity) $n$ if $\parder{^k f}{y^k}(y_0)=0$, $\forall k \leq n$.
\end{defi}

For the proof of the following properties we refer to the classical result by Hartman and Wintner \cite{HW} and to a recent improvement in \cite{HHOT} (Theorem 2.1).
\begin{teo}\label{theorem_nodal_set_f}
Let $f$ be a non trivial solution of \eqref{eq:elliptic_generic}.
\begin{itemize}
\item[(i)] The interior singular points of $f$ are isolated and have finite multiplicity $n \in \mathbb{N}$ ($n \geq 1$).
\item[(ii)] The nodal set of $f$ is the union of finitely many connected arcs which, for a suitable choice of the parametrization, are locally $C^{1,\a}$ for every $\a\in(0,1)$. Moreover such arcs have endpoints either at $\partial \Omega$ or at interior singular points.
\item[(iii)] If $f$ has a zero of order $n$ at the origin, then there exists a function $\tilde\xi \in C^{0,\alpha}(D, \C), \forall \alpha\in (0,1)$, with $\tilde\xi(0)=0$, such that
\begin{equation*}
f(\rho,\phi)=\frac{\rho^{n+1}}{n+1}\Big\{ c_{n+1} \cos[(n+1)\phi] + d_{n+1} \sin[(n+1)\phi]+\tilde\xi(\rho,\phi) \Big\},
\end{equation*}
where $y=\rho e^{i \phi}$. Equivalently, there exists $\xi \in C^{0,\alpha}(D, \C), \forall \alpha\in (0,1)$, such that
\[
2\parder{f}{y}(y)=y^n\xi(y), \qquad \xi(0)=c_{n+1}-i d_{n+1}.
\]
\item[(iv)] If $f$ has a zero of order $n$ at the origin, then for every $k \leq n$ the following Cauchy formula is available
\[
\frac{2}{y^k}\parder{f}{y}(y)=-\frac{i}{\pi}\int_{\partial D} \frac{1}{z^k(z-y)}\parder{f}{z}(z) \ dz
+\frac{1}{2\pi}\int_D \frac{-\Delta f(z)}{z^k(z-y)} \ dz_1dz_2
\]
where the first integral is a complex line integral, whereas the second one is a double integral in the real variables $z_1,z_2$. In particular, an expression for the first non zero coefficients of the expansion of $f$ at the origin is
\begin{equation*}
c_{n+1}-i d_{n+1} =2i\int_{\partial\Omega} \parder{f}{y} \ h_{n+1}  \ dy
-\int_\Omega -\Delta f \ h_{n+1} \ dy_1dy_2
\end{equation*}
where
\begin{equation*}
h_{n+1}(y)=-\frac{1}{2\pi}\frac{1}{y^{n+1}}.
\end{equation*}
\end{itemize}
\end{teo}

We shall now prove that, under the assumptions we are considering, we can still recover similar properties for the magnetic Schr\"odinger equation. The local behavior of the nodal lines is unaltered far from the singular point $a$, but the global behavior undergoes meaningful changes.

\begin{proof}[Proof of Theorem \ref{theorem_nodal_set_U_1}]
Choose $\mpAa=\nabla\Thetaa$ as in Lemma \ref{lemma_potenziale2}, so that all the results of Subsection \ref{subsec_related_real_elliptic_problems} hold true. Let us consider in particular the function $\uuno:\Omegaa\to\R$ as in Lemma \ref{lemma_definition_uuno}. Being solution of a real elliptic equation in a bounded domain, it clearly satisfies the properties collected in the previous theorem. Now, it suffices to notice that $\Pi_y$ is locally holomorphic in every open set which does not contain the point $(a,0)\in\Sigmaa$, hence the local properties of the nodal lines are preserved in the composition and $U$ satisfies Theorem \ref{theorem_nodal_set_f} at every singular point different from $a$. In order to prove that there is at least one nodal line ending at $a$, observe that $\Omegaa$ is symmetric with respect to the origin and $\uuno$ is odd. This implies that the nodal set of $\uuno$ is also symmetric, in particular there are at least two arcs of nodal line having an endpoint at the origin.

In order to prove (iii), let us consider the asymptotic expansion of $\uuno$ near the origin\footnote{The point $a\in D$ corresponds to the origin in $\Omegaa$.}. Since $\uuno$ is odd, Theorem \ref{theorem_nodal_set_f} (iii) gives, for some odd $k \geq 1$
\begin{equation*}
\uuno (\rho,\phi)=\frac{\rho^{k}}{k}\big[ c_{k} \cos(k\phi) + d_{k} \sin(k\phi) \big] + o(\rho^{k}),
\end{equation*}
where we used the notation of Definition \ref{def_polar_coordinates}. From the definition of $\uuno$ we can recover an expression\footnote{Here we wrote $u(x)$ instead of $u(x,y)$ since $\Sigmaa$ is endowed with the locally flat metric induced by $\Pi_x$.} for $u$:
$$
u(r_a,\theta_a)=\frac{{r_a}^\frac{k}{2}}{k}\left[ c_{k} \cos\left(\frac{k}{2}\theta_a\right) + d_{k} \sin\left(\frac{k}{2}\theta_a\right) \right] + o({r_a}^\frac{k}{2}).
$$
This last expression is well defined on $\Sigmaa$, since the complex square root function is continuous on the twofold covering manifold. Finally, \eqref{U_u} provides the corresponding expression for $\tilde{U}$ which, being symmetric, can be projected on $D$, providing \eqref{sviluppo_U}.
\end{proof}

\begin{rem}
It comes from the previous proof that the complex phase $\Thetaa$ which appears in (\ref{sviluppo_U}) is precisely the function $\Thetaa$ defined in Lemma \ref{lemma_potenziale2}.
\end{rem}

\begin{prop} \label{prop_nodal_set_U}
The first non zero coefficients of the asymptotic formula (\ref{sviluppo_U}) can be expressed as
\begin{eqnarray} \label{coefficienti_U}
c_{k}-i d_{k} = 4i \int_{\partial D} G_{k} \left(\parder{U}{x}-i \mpAa U \right)  \ dx 
-2 \int_D G_{k} (i\nabla + \mpAa)^2 U  \ dx_1dx_2,
\end{eqnarray}
where the first integral is a complex line integral, whereas the second one is a double integral in the real variables $x_1,x_2$, and
\begin{equation} \label{G}
G_{k}=-\frac{1}{2\pi}\frac{e^{-i\Thetaa}}{(x-a)^\frac{k}{2}}.
\end{equation}
\end{prop}

\begin{proof} 
Such as in the previous proof, let us consider the function $\uuno$: Theorem \ref{theorem_nodal_set_f}, (iv) gives
$$
c_{k}-i d_{k} =2i\int_{\partial\Omegaa} \parder{\uuno}{y} \ h_{k}  \ dy
-\int_{\Omegaa} -\Delta \uuno \ h_{k} \ dy_1dy_2,
\quad \text{ with } \ h_k(y)=-\frac{1}{2\pi}\frac{1}{y^k}
$$
Let us remark again that the first integral is a complex line integral, whereas the second one is a double integral in real variables, hence by conformal invariance we obtain
$$
c_{k}-i d_{k} = 2i\int_{\partial\Sigmaa} \parder{u}{x} \ g_{k}  \ dx
-\int_{\Sigmaa} -\Delta u \ g_{k} \ dx_1dx_2,
$$
where $g_{k}(x,y)=h_{k} \circ \Pi_y (x,y)$. By taking the complex derivative in \eqref{U_u} we obtain
\[
\parder{u}{x}=e^{-i\Thetaa}\left(\parder{}{x}-i \parder{\Thetaa}{x} \right)\tilde{U},
\]
and
\begin{eqnarray*}
-\Delta u &=& -4 \parder{}{\bar{x}}\parder{u}{x}
=-4 e^{-i\Thetaa} \left(\parder{}{\bar{x}}-i\parder{\Thetaa}{\bar{x}}\right)
\left(\parder{}{x}-i\parder{\Thetaa}{x}\right) \tilde{U} \\
&=&4 e^{-i\Thetaa}\left[-\parder{}{\bar{x}}\parder{}{x}+i\left(\parder{\Thetaa}{x}\parder{}{\bar{x}}
+\parder{\Thetaa}{\bar{x}}\parder{}{x}\right) +\parder{\Thetaa}{\bar{x}}\parder{\Thetaa}{x}\right] \tilde{U} \\
&=& e^{-i\Thetaa} (i \nabla +\tilde{\mpAa})^2 \tilde{U}.
\end{eqnarray*}
By replacing the last expressions in the integral, we obtain
$$
c_{k}-i d_{k}=2i\int_{\partial\Sigmaa} e^{-i\Thetaa} g_{k} \left(\parder{\tilde{U}}{x}-i \tilde{U} \parder{\Thetaa}{x} \right) \ dx 
- \int_{\Sigmaa} e^{-i\Thetaa} g_{k} (i\nabla + \tilde{\mpAa})^2 \tilde{U} \ dx_1dx_2.
$$
In order to conclude the proof it is sufficient to define
$$
G_{k}(x,y)=e^{-i\Thetaa}g_{k}(x,y)=e^{-i\Thetaa}h_{k} \circ \Pi_y (x,y),
$$
and then to observe that both integrands are symmetric on $\Sigmaa$, therefore the last expression can be projected on $D$.
\end{proof}

Let us now turn to the particular case $\Gamma\in\mathcal G$.

\begin{proof}[Proof of Theorem \ref{theorem_nodal_set_U_2}]
We claim that the nodal arcs of $U$ can not be closed curves. In fact $\uuno$ satisfies the maximum principle, since assumption \eqref{hypothesis_V} is preserved by conformal transformations. Thus the nodal lines of $\uuno$ can not be closed curves (by the unique continuation property for real elliptic equations), and this property is preserved by the projections $\Pi_x, \Pi_y$. Now, since $\Gamma$ vanishes exactly three times on $\partial D$, then by simple geometric considerations we infer that there can be at most three nodal lines. The second part of the statement can be obtained similarly, by analyzing the behavior of the nodal set of $\uuno$; the number of nodal lines depends in particular on the form of the boundary trace. Since $\Gamma\in\mathcal{G}$ then, by Lemma \ref{lemma_psi_square}, $\gamma=\psi^{-1/2}\Gammatilde$ is real valued and antisymmetric. Hence in particular $|\psi^{-1/2}\Gammatilde|$ can be projected on $D$ and, by definition of $\mathcal G$, there exists $(\gamma_1,\gamma_2,\gamma_3)\in g$ such that $|\psi^{-1/2}\Gammatilde|=|\Gamma|=\sum_{i=1}^3 \gamma_i$. If $\Phi$ is the global chart of $\Sigmaa$ introduced in Proposition \ref{prop_sigma}, then $\gamma$ has the form
\begin{eqnarray*}
\gamma\circ\Phi(\theta)=\sum_{i=1}^3 \sigma_i \gamma_i \quad \text{for }\theta\in [0,2\pi) \qquad
\gamma\circ\Phi(\theta)=-\gamma\circ\Phi(\theta-2\pi) \quad \text{for }\theta\in [2\pi,4\pi)
\end{eqnarray*}
for some $\sigma_i=\pm1$. It has just been shown that the nodal regions are connected, hence, depending on the signs $\sigma_i$, one obtains the configurations described in the statement.
\end{proof}

\section{Criticality and local uniqueness}\label{section_local_uniqueness}
This section contains the technical part of the paper. We shall prove Theorems \ref{theorem_eigenvalues}, \ref{theorem_criticality} and a local uniqueness result (a local version of Theorem \ref{theorem_global_uniqueness}). Throughout this section $\Gamma \in \mathcal{G}$ is fixed and $\mpAa=\nabla\Thetaa$ is the potential chosen in Lemma \ref{lemma_potenziale2}. Since the position of the singularity is now a variable of the problem, we will denote by $\Ua$ the corresponding solution of (\ref{cauchy_pb_A}), stressing the dependence on the parameter $a$. Similarly, $\ua$ will be the function defined in Proposition \ref{prop_U_u} (with boundary trace $\gamma$ as in the previous proof). Moreover we will denote by $\uunoa$ and $\uduea$ the functions introduced in Lemma \ref{lemma_definition_uuno} and \ref{lemma_definition_udue} respectively.

\subsection{Preliminary estimates}

\begin{lemma}\label{lemma_C^1_estimates}
Let $\Omegatilde$ be a compact subset of $D$. For every $\a\in(0,1)$ there exists a constant $C>0$ such that
$$
\|\uduea-\udueO \|_{C^1(\Omegatilde)} \leq C |a|^\a, \qquad \forall \, a\in D.
$$
The same estimate holds for $\uunoa$, whenever $\Omegatilde\subset \Omega_{a} \cap \Omega_0$.
\end{lemma}

\begin{proof}
By standard imbeddings, $V \in C^{0,\a}(D), \,\forall \,\alpha \in (0,1)$ and $\gamma \in C^{0,1}(\partial D)$. By using Remark \ref{remark_Vdue} and remembering that the M\"obius transformation $T_a(y^2)$ is regular, it is easy to see that for every $\a\in(0,1)$ there exists a constant $C>0$ such that
$$
\|\Vdue_{a} -\Vdue_{0} \|_{L^\infty(D)} + \|\gammadue_{a} -\gammadue_{0} \|_{L^\infty(\partial D)}\leq C|a|^\a, \quad \forall \, a\in D.
$$
On the other hand, $\uduea$ is solution of an elliptic problem, hence by standard regularity results it holds
$$
\|\udue_{a}-\udue_{0}\|_{W^{2,2}(\Omegatilde)} \leq C|a|^\a,
$$
which gives the first part of the statement. Finally observe that $\uunoa$ is the composition of $\uduea$ with a regular function (by Lemma \ref{lemma_definition_udue} again), hence the same estimate holds, whenever it is well defined.
\end{proof}

In case of triple point configuration, the previous estimates can be improved as follows.

\begin{lemma} \label{lemma_L^infty_estimates}
Suppose that the origin is a triple point for $\Gamma$. There exists $\epsilon>0$ such that for every $\alpha\in (1/2,1)$ there exists $C>0$ such that it holds
\begin{equation*}
\| \ua \|_{L^\infty(D_{2|a|})} \leq C |a|^\a, \qquad
\|\nabla \ua \|_{L^\infty(\partial D_{2|a|})} \leq C |a|^{\a-1/2},
\end{equation*}
for every $a\in D$ with $|a|<\epsilon$.
\end{lemma}

\begin{proof}
Since the origin is a triple point, the asymptotic expansion (\ref{sviluppo_U}) of $\UO$ holds with $k=3$, hence $|\uuno_0|\leq Cr^3$, $|\nabla\uuno_0|\leq C' r^2$ in $D_r$. Let $\epsilon>0$ be such that the disk of radius $|3a|^{1/2}$ is contained in $\Omegaa$, whenever $|a|<\epsilon$. Then Lemma \ref{lemma_C^1_estimates} implies the existence of $C>0$ such that
\begin{equation*}
\|\uunoa\|_{L^\infty} + \|\nabla \uunoa\|_{L^\infty} \leq C |a|^\a \quad \textnormal{in } D_{\sqrt{3|a|}},
\end{equation*}
whenever $|a|<\epsilon$. Since $|x|\leq2|a|$ implies $|y| \leq |3 a|^{1/2}$ in $\Sigmaa$, this immediately gives the first inequality of the statement. Now, by definition it holds
\[
|\nabla \ua(y^2+a,y)|=\frac{|\nabla \uunoa(y)|}{2|y|}.
\]
Observe that $|x|=2|a|$ implies $|a|^{1/2}\leq |y| \leq |3 a|^{1/2}$ in $\Sigmaa$, hence we finally obtain
\[
\| \nabla \ua\|_{L^\infty(\partial D_{2|a|})} \leq C\frac{|a|^\alpha}{2 |a|^{1/2}},
\]
which is the second inequality.
\end{proof}

\subsection{Proof of Theorem \ref{theorem_criticality}.}\label{subsection_proof_theorem_criticality}
Let us first establish, under regularity assumptions on the boundary trace, the differentiability of the function $\phi(a)$ (which is defined in the statement of Theorem \ref{theorem_criticality}).

\begin{prop}\label{prop_diffentiability_phi}
Assume that $\gamma \in C^{1,1}(\partial \Omega)$, then the function $\varphi(a)$ defined in Theorem \ref{theorem_criticality} is differentiable.
\end{prop}

\begin{proof}
Notice that we can rewrite $\phi(a)$ in the following way:
$$
\phi(a)=\int_{D}\left(\vert (i\nabla+\mpAa)\Ua \vert^2+V\vert \Ua \vert^2\right) dx_1 dx_2
=\frac{1}{2}\int_{D}\left(\vert \nabla \uduea \vert^2+\Vduea (\uduea)^2\right) dy_1 dy_2.
$$
Let us start showing the existence of the partial derivatives of $\phi(a)$; without loss of generality we can consider the derivative in the direction $a=(a,0)$, centered at the origin. By regularity assumptions and Remark \ref{remark_Vdue}, there exist the weak derivatives $\parder{\Vduea}{a}\in L^\infty(D)$ and $\parder{^2\gammaduea}{a^2} \in L^\infty(\partial D)$. If we prove that
\begin{equation} \label{eq:differentiability_phi}
\lim_{a \to 0} \left\| \frac{\Vduea-\VdueO}{a}- {\parder{\Vduea}{a}}{\big\vert_{a=0}}\right\|_{L^2(D)} =
\lim_{a \to 0} \left\| \frac{\gammaduea-\gammadueO}{a}- {\parder{\gammaduea}{a}}{\big\vert_{a=0}}\right\|_{H^1(\partial D)} =0,
\end{equation}
then standard regularity results for elliptic equations ensure the existence of $w \in H^1(D)$, solution of the following equation
\begin{equation*}
\left\{ \begin{array}{ll}
-\Delta w +\VdueO w +{\parder{\Vduea}{a}}{\big\vert_{a=0}} \udue_0 = 0 \quad & \mathrm{in} \ D  \\
w= \parder{\gammaduea}{a}{\big\vert_{a=0}} & \mathrm{on} \ \partial D,
\end{array} \right.
\end{equation*}
such that 
$$
\lim_{a\to0} \left\| \frac{\uduea-\udueO}{a}-w \right\|_{H^1(D)}=0.
$$
This implies the existence of the partial derivative
$$
\parder{\phi}{a}(0)=\int_{D} \left( \nabla\udueO\cdot\nabla w+\VdueO\udueO w 
+\frac{1}{2} {\parder{\Vduea}{a}}{\big\vert_{a=0}} (\udueO)^2 \right) dy_1 dy_2.
$$
Hence let us prove (\ref{eq:differentiability_phi}). In order to simplify notations we denote here $R(a,y):=T_a(y^2)$, where $T_a$ is defined in (\ref{eq:moebius}). It is sufficient to estimate the following quantity (as $a\to0$) since, by Remark \ref{remark_Vdue} again, the other terms are regular
\begin{multline*}
\left\| \frac{V(R(a,y))-V(R(0,y))}{a}- {\parder{V(R(a,y))}{a}}{\big\vert_{a=0}} \right\|_{L^2(D)} \leq \\
\leq \int_0^1 \left\| \nabla_x V(R(ta,y))\parder{R(ta,y)}{a} -\nabla_x V(R(0,y))\parder{R(a,y)}{a}{\big\vert_{a=0}}\right\|_{L^2(D)} d t.
\end{multline*}
By Lusin's theorem, the integrand converges to zero outside an arbitrarily small set. Then, by applying Lebesgue convergence theorem, we obtain the first relation in (\ref{eq:differentiability_phi}). The second one can be proved in a similar way, implying the existence of the partial derivatives of $\phi(a)$. In order to prove differentiability we test the equation for $\udueO$ with $w$, obtaining
$$
\parder{\phi}{a}(0)=-\int_{\partial D} w \nabla\udueO \cdot \nu d\sigma + \frac{1}{2}
\int_{D}{\parder{\Vduea}{a}}{\big\vert_{a=0}} (\udueO)^2 dy_1 dy_2.
$$
The continuity of this function, with respect to $a$, comes from Lemma \ref{lemma_C^1_estimates} and from the regularity of $\gamma$.
\end{proof}

The proof of Theorem \ref{theorem_criticality} will be divided into two steps. First we need to show that the triple points are critical for the function $\phi(a)$; it will be an immediate consequence of the following proposition.

\begin{prop}\label{prop_triple_point_is_critical}
Suppose that the origin is a triple point for $\Gamma$, then
\[
\lim_{|a|\to 0} \; \frac{Q_{\mpAa,V}(\Ua)-Q_{\mpAO,V}(\UO)}{|a|}=0.
\]
In particular, by the coercivity assumption, the following holds
\[
\lim_{|a|\to 0} \; \frac{|| \Ua -\UO ||_{L^{2}(D)}}{|a|} = 0.
\]
\end{prop}

\begin{proof}
We split the energy function into the sum of two integrals:
\begin{eqnarray*}
Q_{\mpAa,V}(\Ua)-Q_{\mpAO,V}(\UO)&=&\int_D \left( |(i\nabla +\mpAa)\Ua |^2 -|(i\nabla +\mpAO)\UO |^2 +V(|\Ua|^2-|\UO|^2) \right)dx_1 dx_2 \\
&=&I + II,
\end{eqnarray*}
where I is the integral in the annulus $D\setminus D_{2|a|}$ and II is the integral in the ball $D_{2|a|}$.

As it concerns the integral in the exterior annulus, the key observation is that both $\uO$ and $\ua$ are well defined in the twofold covering manifold $\Sigma_0\setminus\Pi_x^{-1}(D_{2|a|})$, since the domain $D\setminus D_{2|a|}$ does not contain any singularity. The difference function satisfies the equation
\begin{equation*}
\left\{ \begin{array}{ll}
-\Delta(\ua-\uO)+\tilde{V}(\ua-\uO)=0 \quad & \mathrm{in} \ \Sigma_0\setminus\Pi_x^{-1}(D_{2|a|})  \\
\ua-\uO= 0 & \mathrm{on} \ \partial \Sigma_0,
\end{array} \right.
\end{equation*}
therefore Lemma \ref{lemma_L^infty_estimates} gives, for $a$ sufficiently small,
\begin{eqnarray*}
\int_{D\setminus D_{2|a|}}\left( |\nabla(\ua-\uO)|^2+V (\ua-\uO)^2 \right) dx_1 dx_2 &\leq& 
\int_{\partial D_{2|a|}}|\ua-\uO| \left|\parder{}{\nu}(\ua-\uO)\right| d\sigma \\
&\leq& 4\pi|a| \sup_{\partial D_{2|a|}}\{ |\nabla(\ua-\uO)| |\ua-\uO| \} \\
&\leq& C |a|^{1/2+2\a}.
\end{eqnarray*}
By choosing $\a=7/8$ and using the coercivity assumption (\ref{hypothesis_V}) we infer
\[
\int_{D\setminus D_{2|a|}} |\nabla(\ua-\uO)|^2dx_1 dx_2 + \int_{D\setminus D_{2|a|}} (\ua-\uO)^2 dx_1 dx_2 \leq C|a|^{9/4}.
\]
On the other hand by conformal invariance it holds
\[
I\leq C \left(\int_{D\setminus D_{2|a|}}|\nabla(\ua-\uO)|^2 dx_1 dx_2 \right)^{1/2} +
C' \left(\int_{D\setminus D_{2|a|}}(\ua-\uO)^2 dx_1 dx_2 \right)^{1/2},
\]
which, together with the previous inequality, gives $I \leq C|a|^{9/8}$.

As it concerns the integral in $D_{2|a|}$, we proceed in a similar way and apply Lemma \ref{lemma_C^1_estimates}, obtaining
\begin{eqnarray*}
II &\leq&  C \left(\int_{D_{\sqrt{2|a|}}}|\nabla(\uduea-\udueO)|^2 dx_1 dx_2\right)^{1/2} +
C' \left(\int_{D_{\sqrt{2|a|}}}(\uduea-\udueO)^2 dx_1 dx_2\right)^{1/2}\\
&\leq& |a|^{1/2+\a}.
\end{eqnarray*}
By choosing $\a=5/8$ here, we finally obtain
$$
\lim_{|a|\to 0} \frac{I+II}{|a|}=\lim_{|a|\to 0} |a|^{1/8} =0,
$$
which concludes the proof.
\end{proof}

The next result deals with the case there is no multiple clustering and shows that the converse of the previous proposition also holds.

\begin{prop}\label{prop_energy_decrease}
Suppose that the origin is not a triple point for $\Gamma$. There exist $\bar h, C>0$ such that for every $0<h<\bar h$ there exist $a\in D$ with $|a|=h$ and a function $\Za \in \HAD$ such that
\[
\lim_{|a|\to 0} \; \frac{Q_{\mpAa,V}(\Za)-Q_{\mpAO,V}(\UO)}{|a|} < -C \qquad\text{and}\qquad
\lim_{|a|\to 0} \; \frac{\|\Za\|_{L^2(D, \C)}^2-\|\UO\|_{L^2(D,\C)}^2}{|a|}=0.
\]
\end{prop}

\begin{proof}
By Theorems \ref{theorem_nodal_set_U_1}, \ref{theorem_nodal_set_U_2}, $\UO$ has the following asymptotic expansion around the origin\footnote{Here, respect to equation (\ref{sviluppo_U}), we have set $\alpha=\arctan(d_1/c_1)$ and $C=c_1/ \cos\alpha \neq 0$.}
\begin{equation}\label{critical_points_phi_sviluppo_asintotico}
\UO(r,\theta)=C e^{i\ThetaO} r^{\frac{1}{2}} \cos(\frac{\theta}{2}-\alpha) + o(r^{\frac{1}{2}}),
\end{equation}
for some $C\neq0$. In particular there is exactly one nodal line ending at the origin and there exists $\bar h$ such that the disk $D_{\bar h}$ does not intersect any other nodal line of $\UO$. For every $h<\bar h$ let $w:D_h \to \R$ be the (nonnegative) solution of
\begin{eqnarray*}
\left\{ \begin{array}{ll}
        -\Delta w +V w=0 & \text{ in } D_h \\
	    w=|\uO| & \text{ on } \partial D_h.
        \end{array}
\right.
\end{eqnarray*}
Let now $a \in \mathcal N (\UO)$ be the unique point with $|a|=h$. We define a new function $\za: \Sigmaa\to \R$ as
\begin{eqnarray*}
\za(x,y)=\left\{ \begin{array}{ll}
\sigma(x,y) \tilde{w}(x,y) & (x,y) \in \Pi_x^{-1}(D_h) \\
\sigma(x,y) |\tilde{u}_0| & (x,y) \in \Sigmaa \setminus \Pi_x^{-1}(D_h),
\end{array}
\right.
\end{eqnarray*}
where $\sigma(x,y)=\pm 1$ in such a way that $\za$ is antisymmetric on $\Sigmaa$. If $\Thetaa: \Sigmaa\to \R$  is defined as in 
Lemma \ref{lemma_potenziale2}, then $\Za=\Pi_x(e^{-i\Thetaa}\za) \in \HAD$ and
\[
Q_{\mpAa,V}(\Za)-Q_{\mpAO,V}(\UO)= \int_{D_h} \left[|\nabla w|^2-|\nabla \uO|^2+V(w^2-\uO^2)\right] dx_1 dx_2.
\]
In order to estimate the limit as $|a|=h\to0$, we perform the following change of variables
$$
\uh(y)= \frac{1}{\sqrt{h}}\uO(h y^2, \sqrt{h}y), \quad
\wh(y)= \frac{1}{\sqrt{h}} w(h y^2).
$$
These functions satisfy the rescaled problems
\begin{eqnarray*}
\left\{ \begin{array}{ll}
        -\Delta \uh + h V(h y^2) \uh =0 & \text{ in } D \\
	\uh=\frac{1}{\sqrt{h}}\uO(h y^2, \sqrt{h}y) & \text{ on } \partial D,
        \end{array}
\right.
\quad
\left\{ \begin{array}{ll}
        -\Delta \wh + h V(h y^2) \wh =0 & \text{ in } D \\
	\wh=|\uh | & \text{ on } \partial D.
        \end{array}
\right.
\end{eqnarray*}
Moreover by \eqref{critical_points_phi_sviluppo_asintotico}, $\uh$ satisfies the asymptotic expansion
$$
\uh(\rho,\phi)=C \rho \cos(\phi-\alpha) + o(\sqrt{h}\rho),
$$
where as usual $y=\rho e^{i\phi}$. This ensures the existence of a limit function $\uinfty$ such that
\begin{eqnarray*}
\left\{ \begin{array}{ll}
        -\Delta \uinfty =0 & \text{ in } D \\
	\uinfty(\rho, \phi)= C \cos(\phi-\alpha) & \text{ on } \partial D
        \end{array}
\right.
\ \text{ and } \
\| \uh - \uinfty\|_{C^1(\overline{D})} \to 0, \text{ as } h \to 0.
\end{eqnarray*}
As a consequence, $\| |\uh| - |\uinfty| \|_{W^{1,p}(\partial D)} \to 0, \ \forall p \in (1,+\infty)$, which implies
\begin{eqnarray*}
\left\{ \begin{array}{ll}
        -\Delta \winfty =0 & \text{ in } D \\
	\winfty(\rho, \phi)= |C \cos(\phi-\alpha)| & \text{ on } \partial D
        \end{array}
\right.
\ \text{ and } \
\| \wh - \winfty\|_{H^1(D)} \to 0, \text{ as } h \to 0.
\end{eqnarray*}
Therefore we have obtained
\begin{eqnarray*}
\lim_{|a|\to0}\frac{Q_{\mpAa,V}(\Za)-Q_{\mpAO,V}(\UO)}{|a|} &=&
\lim_{h\to 0}\frac{1}{h}\int_{D_h} \left[|\nabla w|^2-|\nabla \uO|^2+V(w^2-\uO^2)\right] dx_1 dx_2\\
&=& \int_{D} \left[|\nabla \winfty |^2-|\nabla \uinfty|^2 \right] dy_1 dy_2.
\end{eqnarray*}
The second inequality in the statement is now proved, let us go to the first one. By choosing the coordinates in such a way that $\alpha=0$ in \eqref{critical_points_phi_sviluppo_asintotico}, we have
$$
\uinfty=C \rho \cos(\phi), \quad \int_{D}|\nabla \uinfty|^2 dy_1 dy_2 = C^2 \pi,
$$
and
$$
\winfty=\frac{2|C|}{\pi} \left[ 1+ 2\sum_{n=1}^\infty \frac{(-1)^n}{1-4 n^2} \rho^{2n} \cos(2n\phi) \right],
$$
which gives
\[
\int_{D}|\nabla \winfty|^2 dy_1 dy_2 =
\frac{32 C^2}{\pi}\sum_{n=1}^\infty \frac{n}{(1-4n^2)^2} \leq \frac{44 C^2}{9 \pi},
\]
concluding the proof (since $C\neq 0$).
\end{proof}

\begin{proof}[Proof of Theorem \ref{theorem_criticality}]
Assume first the origin is a triple point for $\Gamma$ (this is not restrictive since we can always apply the the conformal map $T_a$ defined in (\ref{eq:moebius})). We have shown in Proposition \ref{prop_triple_point_is_critical} that
$$
\phi'(0)=\lim_{|a|\to 0} \frac{Q_{\mpAa,V}(\Ua)-Q_{\mpAO,V}(\UO)}{|a|}=0,
$$
hence the origin is a critical point of $\phi(a)$. On the other hand if the origin is not a triple point for $\Gamma$, then by definition $\phi(a) \leq Q_{\mpAa,V}(\Za)$, where $a$ is sufficiently small and $\Za$ is defined in Proposition \ref{prop_energy_decrease}. Therefore there exists $C>0$ such that
\[
\phi'(0) \leq \lim_{|a|\to 0} \frac{Q_{\mpAa,V}(\Za)-Q_{\mpAO,V}(\UO)}{|a|} <-C
\]
and in this case the origin is not a critical point of $\phi(a)$.
\end{proof}

\subsection{Proof of Theorem \ref{theorem_eigenvalues}}
Let us fix the notations as follows. Let $\lambda_a^k$ be the $k$--eigenvalue of the operator $H_{\mpAa,0}$ with Dirichlet boundary conditions in $D$ ($k$ positive integer). Correspondingly, we denote by $\psi_a^k$ any associated eigenfunction normalized in the $L^2$--norm.
\begin{rem}\label{rem_nodal_lines_eigenfunctions}
There exists a basis of $K$--real eigenfunctions of $H_{\mpAa,0}$ with Dirichlet boundary conditions. Moreover every eigenfunction, being a complex multiple of a $K$--real function, satisfies properties (i)-(ii)-(iii) in Theorem \ref{theorem_nodal_set_U_1}.
\end{rem}
We recall that
\[
\lambda_a^k =\inf\left\{ \frac{\int_D |(i\nabla+\mpAa)U|^2 dx_1 dx_2}{\int_D |U|^2 dx_1 dx_2} : \ U\in H^1_0(D), \int_D U \overline \Psi^j_a dx_1 dx_2 =0, \ \forall \ 1\leq j\leq k-1 \right\}.
\]

The proof of Theorem \ref{theorem_eigenvalues} is similar to the one of Theorem \ref{theorem_criticality}, except for the fact that the function $\Za$ introduced therein needs to be orthogonal to the first $k-1$ eigenfunctions. This is precisely the content of the following proposition.

\begin{prop}\label{prop_energy_decrease_orthogonal}
Assume that the nodal set of $\Psi_0^k$ does not present multiple clustering at the origin. There exist $\bar h, C>0$ such that for every $0<h<\bar h$ there exist $a\in D$ satisfying $|a|=h$ and a function $\Za \in \HAD$ such that
\[
\lim_{|a|\to 0} \; \frac{Q_{\mpAa,0}(\Za)-Q_{\mpAO,0}(\Psi^k_0)}{|a|} < -C, \qquad
\lim_{|a|\to 0} \; \frac{\|\Za\|_{L^2(D,\C)}^2-\|\Psi^k_0\|_{L^2(D,\C)}^2}{|a|}=0
\]
and moreover, if $k\geq 2$, it holds $\int_D \Za \overline \Psi_a^j=0$ for every $1\leq j\leq k-1$.
\end{prop}

\begin{proof}
Let us start assuming that the first $k$ eigenvalues of $H_{\mpAa,0}$ are simple. We shall perform the proof for $k=1,2$, it should be clear how to modify it in the remaining cases. Let us choose $\epsilon>0$ so small that for every $\alpha \in [-\epsilon,\epsilon]$ the linear combination
\[
\Psi(\alpha)=\alpha\Psi^1_0+\sqrt{1-\alpha^2}\Psi^2_0
\]
admits a local expansion
\[
\Psi(\alpha) =e^{i\frac{\theta}{2}} r^{\frac{1}{2}}\left[ c_1(\alpha) \cos\frac{\theta}{2} + d_1(\alpha) \sin\frac{\theta}{2} \right] + o(r^{\frac{1}{2}}),
\]
where $c_1(\alpha),d_1(\alpha)$ are real constant such that $c_1(\alpha)^2+d_1(\alpha)^2\geq\delta$ for some $\delta>0$. Such an $\epsilon$ clearly exists by Remark \ref{rem_nodal_lines_eigenfunctions} and by assumption. Therefore there exists $\bar h>0$ such that the disk $D_{\bar h}$ intersects exactly one nodal line of $\Psi(\alpha)$ for every $\alpha \in [-\epsilon,\epsilon]$. For $0\leq h<\bar h$ let $a(\alpha, h)$ be the only point satisfying $a(\alpha,h) \in \mathcal N(\Psi(\alpha))$ and $|a(\alpha,h)|=h$; clearly $a(\alpha,0)=0$. Let now
\begin{eqnarray*}
\left\{ \begin{array}{ll}
        -\Delta w = \lambda^2_0 w& \text{ in } D_h \\
	    w=|\Psi(\alpha)| & \text{ on } \partial D_h.
        \end{array}
\right.
\end{eqnarray*}
and define $\Za=Z_{a(\alpha,h)}$ similarly to Proposition \ref{prop_energy_decrease}. Proceeding as therein, we immediately obtain the existence of a constant $C>0$ such that
\[
\lim_{|a|\to 0} \; \frac{Q_{\mpAa,0}(\Za)-Q_{\mpAO,0}(\Psi(\alpha))}{|a|} < -C, \qquad
\lim_{|a|\to 0} \; \frac{\|\Za\|_{L^2(D)}^2-\|\Psi(\alpha)\|_{L^2(D)}^2}{|a|}=0
\]
for every $\alpha \in [-\epsilon,\epsilon]$ and $0<h<\bar h$. Now, since
\[
\int_D|(i\nabla+\mpAO)\Psi(\alpha)|^2dx_1 dx_2=\alpha\lambda^1_0+\sqrt{1-\alpha^2}\lambda^2_0\leq \lambda^2_0, \qquad \text{and} \qquad \int_D|\Psi(\alpha)|^2dx_1 dx_2=1,
\]
we immediately obtain that the functions $Z_{a(\alpha,h)}$ satisfy the first part of the statement for every $\alpha \in [-\epsilon,\epsilon]$, hence the proof is complete in the case $k=1$. If $k=2$, let us show that we can select $\alpha$ in such a way that the orthogonality condition is satisfied. To this aim consider, for every fixed $0\leq h<\bar h$, the following map
\begin{eqnarray*}
[-\epsilon,\epsilon] &\to& \R\\
\alpha &\mapsto& \int_D Z_{a(\alpha,h)} \overline \Psi^1_{a(\alpha,h)}dx_1 dx_2.
\end{eqnarray*}
It coincides with the identity for $h=0$ and, being the eigenvalues simple, depends continuously on $h$. Therefore, thanks to the Theorem of Borsuk--Ulam, it admits a zero for every $h$ sufficiently small, proving the proposition in case of simple eigenvalues.

In the general case it is possible to perform an arbitrarily small perturbation of the domain in such a way that the eigenvalues are simple on the new domain (see \cite{M}). Moreover, the new domain can be chosen diffeomorphic to disk, hence by conformal invariance this is equivalent to working with generalized eigenfunctions, which minimize the generalized Rayleigh quotient
\[
\frac{\int_D |(i\nabla+\mpAa)U|^2dx_1 dx_2}{\int_D V(x)|U|^2dx_1 dx_2},
\]
where $V(x)>0$ is the Jacobian of the conformal transformation. It only remains to check that the constant $C$ which appears in the statement is independent on the choice of the potential $V(x)$. But this was already proved in Proposition \ref{prop_energy_decrease}, where the constant $C$ is explicitly determined. Finally, by regularity, the energy of the perturbed problem converges to the one of the original problem.
\end{proof}

\begin{proof}[Proof of Theorem \ref{theorem_eigenvalues}]
Assume first that $\Psi^2_0$ does not present multiple intersection of the nodal lines at the origin; let us prove that the origin is not a critical point of the function $a \mapsto \lambda_a^2$. Let $\Za$ be the function introduced in Proposition \ref{prop_energy_decrease_orthogonal}, then
\[
\lambda^2_a \leq \frac{\int_D |(i\nabla+\mpAa)\Za|^2dx_1 dx_2}{\int_D|\Za|^2dx_1 dx_2},
\]
hence a simple manipulation gives (by recalling that $\int_D|\Psi_0^2|^2dx_1 dx_2=1$)
\[
\begin{split}
\lim_{|a|\to 0} \frac{\lambda^2_a -\lambda^2_0}{|a|} \leq 
\lim_{|a|\to 0} \int_D \left(|(i\nabla+\mpAO)\Za|^2 \right.&-\left. |(i\nabla+\mpAO)\Psi^2_0|^2 \right)dx_1 dx_2 \\
&+ \lambda^2_0 \lim_{|a|\to 0} \int_D\left(|\Psi^2_0|^2-|\Za|^2 \right)dx_1 dx_2,
\end{split}
\]
which is strictly negative by virtue of Proposition \ref{prop_energy_decrease_orthogonal}.
\end{proof} 
\begin{proof}[Proof of Corollary \ref{coro_eigenvalues}]
As it concerns the first part of the statement, the only non trivial fact is the differentiability of the eigenvalue as long as it remains simple. This can be seen, for example, by using Riemann mapping theorem composed with a M\"obius map in order  to transform the punctured domain in the standard punctured disk with the origin removed. Next we can square the independent variable and lift the eigenvalues problem to the disk. By conformal invariance, the new eigenvalues problem depends smoothly on the parameter $a$.

Let us turn to the second part of the statement, concerning the first eigenvalue. It is not difficult to show (see \cite{NT}) that
\[
\lim_{a\to\partial \Omega}\lambda_a^1=\lambda_1(\Omega),
\]
whereas the diamagnetic inequality gives $\lambda_a^1>\lambda_1(\Omega)$, hence the map $a\mapsto \lambda_a^k$ must have an interior maximum. In order to conclude the proof, assume by contradiction that the first eigenvalue is simple at this maximum. Then $\lambda_a^k$ would be differentiable with respect to $a$ and Theorem \ref{theorem_eigenvalues} would imply the existence of a multiple junction at $a$. But this is a contradiction, since the first eigenfunction $\Psi^1_a$ has exactly one nodal line. On the other hand, it is proved in \cite{HHO1} that the multiplicity of the first eigenvalue is at most two, which concludes the proof.
\end{proof}

\subsection{Local uniqueness of the triple point.}\label{subsection_local_uniqueness}
Following \cite{ctv5}, we shall base the proof of Theorem \ref{theorem_global_uniqueness} on a local uniqueness result.

\begin{teo} \textnormal{\textbf{(Local uniqueness)}} \label{theorem_local_uniqueness}
Suppose that $\Gamma \in \mathcal{G}$ admits a triple point $a_\Gamma\in D$. Then there exist $\varepsilon > 0, C>0$ such that, for every boundary data $\Lambda \in \mathcal{G}$ satisfying $||\Gamma-\Lambda||_{L^\infty(\partial D)} < \varepsilon$, there exists exactly one $a_\Lambda$ (triple point for $\Lambda$) satisfying $|a_\Gamma -a_\Lambda|< C ||\Gamma-\Lambda||_{L^\infty(\partial D)}$.
\end{teo}

\begin{proof}
We can assume without loss of generality that $a_\Gamma=0$ (by applying the conformal map $T_{a_\Gamma}$ defined in (\ref{eq:moebius})), in such a way that the function $U_0$ (with boundary trace $\Gamma$) has a triple point at the origin. Let now $\Lambda \in \mathcal{G}$ be any boundary trace and let $\Ua$ be the corresponding solution of (\ref{cauchy_pb_A}). Denote by $c_1(a), d_1(a)$ the coefficients of order one in the asymptotic expansion of $\Ua$, as in (\ref{coefficienti_U}). Notice that, if the boundary trace is precisely $\Gamma$ and $a=0$, then, by Theorem \ref{theorem_nodal_set_U_2}, it holds
\begin{equation}\label{c_3^2+d_3^2neq0}
c_1(0)=d_1(0)=0, \qquad c_3(0)^2+d_3(0)^2\neq 0.
\end{equation}
Following \cite{ctv5}, Proposition 3.2, where the authors study the particular case $V=0$, we shall apply the implicit function theorem to the map
\begin{eqnarray} \label{implicit_function}
\mathcal{G} \times D &\to& \R^2 \nonumber \\
(\Lambda, a) &\mapsto& (c_1(a),d_1(a)).
\end{eqnarray}
Theorem \ref{theorem_nodal_set_U_2} again ensures that $a$ is a triple point if and only if $c_1(a)=d_1(a)=0$. Therefore the theorem is proved, provided we can locally solve this equation for $a$ in a neighborhood of $(\Gamma,0)$.

First of all we observe that \eqref{implicit_function} defines a $C^1$ function. Indeed it comes from the proof of Proposition \ref{prop_nodal_set_U} that $(c_1(a),d_1(a))= \nabla_y \uunoa (0)$, and regularity can be proved proceeding as in Proposition \ref{prop_diffentiability_phi}.\footnote{Here we do not need additional regularity on the boundary data, since the estimates are local.} Therefore we only need to show that the $2\times 2$ Jacobian matrix
$$
\nabla_a (c_1(a), d_1(a))\Big\vert_{a=0}
$$
is invertible. By (\ref{coefficienti_U}) it holds
\begin{equation} \label{c1_d1}
c_1(a)-i d_1(a) = 4i\int_{\partial D} \Ga \left(\parder{ }{x}-i \mpAa \right) \Ua \ dx -
2 \int_D \Ga (i\nabla + \mpAa)^2 \Ua  \ dx_1dx_2,
\end{equation}
with
\[
\Ga=-\frac{1}{2\pi}\frac{e^{-i\Thetaa}}{(x-a)^{\frac{1}{2}}}.
\]
Notice that the differential operator commutes with the integral since the functions $\parder{\Ga}{a}(x)\simeq \frac{1}{(x-a)^{3/2}}$ belong to $L^1(D)$ for every $a$. The main difficulty here is that we do not know the behavior of $\Ua$ with respect to the variation of the parameter $a$, therefore we need to manipulate the last expression before differentiating. In order to get rid of the boundary integral in \eqref{c1_d1}, we introduce a new function $\Fa : D \to \C$, solution of the equation
\begin{equation} 
\left\{ \begin{array}{ll}
\left( \parder{}{\xbar} +i \parder{\Thetaa}{\xbar}\right) \Fa =0 \quad & \mathrm{on} \ D \nonumber  \\
\Fa=\Ga & \mathrm{on} \ \partial D. \nonumber
\end{array} \right.
\end{equation}
By applying Green's formula (Lemma \ref{lemma_green}), equation \eqref{c1_d1} becomes
\begin{equation*}
\begin{split}
&c_1(a)-i d_1(a) = 2 \int_D (\Fa-\Ga) (i\nabla + \mpAa)^2 \Ua  \ dx_1dx_2 \\
&= 2 \int_D (\Fa-\Ga) (i\nabla + \mpAO)^2 \UO  \ dx_1dx_2 +2 \int_D (\Fa-\Ga) (V \UO- V \Ua)  \ dx_1dx_2.
\end{split}
\end{equation*}
Now, the last term in the previous equality plays no role in the computation of the derivative, since by Proposition \ref{prop_triple_point_is_critical} it holds
\begin{equation*}
\begin{split}
&2 \left| \int_D (\Fa-\Ga) (V \UO- V \Ua)  \ dx_1dx_2 \right| \\
&\leq 2 ||V||_{L^\infty(D)} ||\Fa-\Ga||_{L^2(D)} ||\UO-\Ua||_{L^2(D)} = o(|a|).
\end{split}
\end{equation*}
Remembering that $\UO$ has a triple point at the origin, we compute
\begin{multline}\label{derivative}
\parder{ }{a}(c_1(a)-i d_1(a))\Big\vert_{a=0}=
\left( \parder{ }{a} +i \parder{\Thetaa}{a} \right) (c_1(a)-i d_1(a)) \Big\vert_{a=0} \\
= \int_D (i\nabla + \mpAO)^2 \UO \left( \left( \parder{ }{a} +i \parder{\Thetaa}{a} \right)\Fa \Big\vert_{a=0} -
\left( \parder{ }{a} +i \parder{\Thetaa}{a} \right)\Ga \Big\vert_{a=0} \right) \ dx_1dx_2.
\end{multline}
We can differentiate $\Ga$ directly since, being integrable, its distributional derivative coincides with the a.e. derivative
$$
\left( \parder{ }{a} +i \parder{\Thetaa}{a} \right)\Ga=-\frac{1}{4\pi}\frac{e^{-i\Thetaa(x)}}{(x-a)^\frac{3}{2}}=
\frac{1}{2} G_{3,a}.
$$
Notice that we obtain a multiple of the function defined in \eqref{G} for $k=3$, which gives information about the asymptotic behavior of the solution at order three. Then we differentiate the equation for $\Fa$
\begin{equation} 
\left\{ \begin{array}{ll}
\left( \parder{}{\xbar} +i \parder{\Thetaa}{\xbar}\right) \left( \parder{ }{a} +i \parder{\Thetaa}{a} \right) \Fa =0 \quad & \mathrm{on} \ D \nonumber  \\
\parder{\Fa}{a}=\parder{\Ga}{a} & \mathrm{on} \ \partial D, \nonumber
\end{array} \right.
\end{equation}
and by using Green's formula again, we recover
\begin{multline*}
\int_D \left( \parder{ }{a} +i \parder{\Thetaa}{a} \right)\Fa \cdot (i\nabla + \mpAO)^2 \UO \ dx_1dx_2 \\
=2i\int_{\partial D} \left( \parder{ }{a} +i \parder{\Thetaa}{a} \right)\Ga \cdot
\left( \parder{ }{x} -i \parder{\Thetaa}{x} \right)\UO \ dx \\
=2i\int_{\partial D} G_{3,a} \cdot \left( \parder{ }{x} -i \parder{\Thetaa}{x} \right)\UO \ dx.
\end{multline*}
Now, by replacing the last expression in \eqref{derivative} we obtain
\begin{equation*}
\parder{ }{a}(c_1(a)-i d_1(a))\Big\vert_{a=0} 
=2i\int_{\partial D} G_{3,0} \left( \parder{ }{x} -i \parder{\ThetaO}{x} \right)\UO \ dx - 
\int_D G_{3,0} (i\nabla + \mpAO)^2 \UO \ dx_1 dx_2,
\end{equation*}
which finally implies
$$
\parder{ }{a}(c_1(a)-i d_1(a))\Big\vert_{a=0}= \frac{1}{2}(c_3(0)-i d_3(0)).
$$
Equation (\ref{c_3^2+d_3^2neq0}) ensures that this quantity does not vanish, hence the implicit function theorem applies to (\ref{implicit_function}) in a neighborhood of $(\Gamma, 0)$ and the theorem is proved.
\end{proof}

\section{Energy minimizing partitions and global uniqueness}\label{section_energy_minimizing_partitions}
In this section we shall prove all the remaining results stated in the introduction, apart from the continuous dependence of the nodal lines with respect to the boundary trace, which will be the object of the last section. All the problems involved are interconnected and strongly related to previous works by Conti, Terracini and Verzini. There is a relation between the triple point configuration of the solutions of (\ref{cauchy_pb_A}), the class $\Sgamma$ introduced in Theorem \ref{theorem_uniqueness_Sgamma} and the optimal partition problem (\ref{minimum_pb}). In order to analyze this relation we shall recall some known results, we refer to \cite{ctv2,ctv3,ctv4,ctv5} for the proofs and for further details. Throughout this section, assumption (\ref{hypothesis_V}) holds.

\begin{defi}
Let $(\gamma_i)\in g$ and correspondingly let $(u_i)\in \Sgamma$. We denote by $\omega_i=\{ u_i>0 \}$. The multiplicity of a point $x \in \overline{D}$ (with respect to $(u_i)$) is
$$
m(x)= \sharp\{ i: \textrm{measure}(\omega_i \cap D_r(x))>0 \quad \forall r>0 \}.
$$
Notice that $1\leq m(x) \leq 3$.
\end{defi}
The properties that we will need can be summarized briethly as follows.
\begin{teo}\label{theorem_properties_Sgamma}
Let $(\gamma_i)\in g$ and correspondingly let $(u_i)\in \Sgamma$, then
\begin{itemize}
\item[(i)] $u_1+u_2+u_3 \in W^{1,\infty}(D)$, hence in particular every $\omega_i$ is open;
\item[(ii)] there exists exactly one point $a \in \overline{D}$ such that $m(a)=3$ with respect to $(u_i)$;
\item[(iii)] let $a$ be as in (ii); if $a \in \partial D$ then $(u_i)$ is the only element of $\Sgamma$;
\item[(iv)] if $m(a)=3$ both with respect to $(u_i)$ and $(v_i)\in \Sgamma$, then $(u_i)=(v_i)$;
\item[(v)] given $a$ as in (ii), consider the twofold covering manifold $\Sigmaa$ and the following functions defined on it
$$
u(x,y):=\sum_{i=1}^3 \sigma(x,y) u_i \circ \Pi_x(x,y), \qquad \gamma(x,y):=\sum_{i=1}^3 \sigma(x,y) \gamma_i \circ \Pi_x(x,y),
$$
where $\sigma(x,y)=\pm 1$ in such a way that they have alternate signs on adjacent supports. Then $u$ satisfies (\ref{cauchy_pb_0}) in $\Sigmaa$ with boundary trace $\gamma$.
\end{itemize}
\end{teo}
With some abuse of notations we will call $a$ the \text{triple point} of $(u_i)$. Properties (ii) and (v) of the previous theorem, together with the analysis performed in Section \ref{sec_gauge}, immediately give
\begin{coro}\label{corollary}
Let $(\gamma_i)\in g$ and let $\Gamma \in \mathcal G$ be such that $|\Gamma|=\sum_{i=1}^3 \gamma_i$. There is a bijection between the elements of $\Sgamma$ having interior triple point (i.e. $a\in D$) and the solutions of (\ref{cauchy_pb_A}) which admit triple point. Moreover, the bijection preserves the nodal set.
\end{coro}
Let us now turn to the relation with the optimal partition problem (\ref{minimum_pb}). It is contained in the following theorem, which is proved in \cite{ctv3}.
\begin{teo}\label{theorem_energy_minimizing_partitions}
Given $(\gamma_i)\in g$, let
\begin{equation*}
\mathcal U=\left\{ (u_1,u_2,u_3) \in (H^1(D))^3: \begin{array}{l}
                                    \ui \geq0\ \textrm{in} \  D, \ \ui=\gamma_i \ \textrm{on} \  \partial D \\
				    \ui\cdot\uj=0 \ \textrm{a.e.} \ x\in D, \ \textrm{for} \ i \neq j
\end{array}
\right\}.
\end{equation*}
The minimization problem
\begin{equation}\label{problema_minimizzazione}
\min_{(u_i)\in \mathcal U}\sum_{i=1}^3\int_D \left( |\nabla u_i|^2 + V u_i^2 \right) dx_1 dx_2
\end{equation}
has a unique solution, which belongs to $\Sgamma$.
\end{teo}


\begin{proof}[Proof of Theorem \ref{theorem_uniqueness_Sgamma} (Sketch)]
First of all by Theorem \ref{theorem_properties_Sgamma}, (iii), we can concentrate on those elements of $\Sgamma$ having interior triple points. Corollary \ref{corollary} together with Theorem \ref{theorem_local_uniqueness} immediately give the following local uniqueness property for $\Sgamma$ (which is the analogous of Proposition 3.2 in \cite{ctv5}):

let $(u_i)\in \Sgamma$ and let $a_\gamma \in D$ be its triple point; there exist $\epsilon, C>0$ such that for every $(\lambda_i) \in g$ satisfying $\sum_{i=1}^3\|\gamma_i-\lambda_i\|_{L^\infty(\partial D)}<\epsilon$, there exists exactly one $a_\lambda \in D$ (triple point for $(\lambda_i)$) such that $|a_\gamma-a_\lambda|< C \sum_{i=1}^3\|\gamma_i-\lambda_i\|_{L^\infty(\partial D)}$.

The fact that $\Sgamma$ consists of exactly one element can now be proved exactly as in \cite{ctv5}, Theorem 1.3; we only sketch the procedure here. Assume by contradiction that $(u_i)\in \Sgamma$ is the solution of (\ref{problema_minimizzazione}) and $(v_i)\neq(u_i) \in \Sgamma$. By Theorem \ref{theorem_properties_Sgamma}, (iv), the triple points of $(u_i)$ and $(v_i)$ can not coincide. The authors apply a blow up procedure to $(v_i)$ in a neighborhood of the triple point and prove convergence. Thus, using the minimality result expressed in Theorem \ref{theorem_energy_minimizing_partitions}, they deduce that $(u_i)$ and $(v_i)$ can be connected by a continuous path of triple point configurations, which contradicts the local uniqueness.

Hence $\Sgamma$ contains exactly one element $(u_1,u_2,u_3)$ which is the solution of (\ref{problema_minimizzazione}). The corresponding supports $\omega_i=\{ u_i>0 \}$, which are open by Theorem \ref{theorem_properties_Sgamma}, (i), are clearly the solutions of (\ref{minimum_pb}).
\end{proof}

\begin{proof}[Proof of Theorems \ref{theorem_minimal_partition} and \ref{theorem_global_uniqueness}] Both the statements are an immediate consequence of Theorem \ref{theorem_uniqueness_Sgamma} and Corollary \ref{corollary}. The second part of Theorem \ref{theorem_global_uniqueness} comes from Theorem \ref{theorem_local_uniqueness}.
\end{proof}

As we mentioned in the introduction, it is possible to give another characterization of these results, in terms of the limiting configuration of the solutions of a competition--diffusion system. In \cite{ctv4} the authors study the following family of elliptic system, with parameter $k$
\begin{equation} \label{sistema_popolazioni}
\left\{ \begin{array}{ll}
-\Delta u_{i,\kappa} + V(x) u_{i,\kappa} =-\kappa u_{i,\kappa} \sum_{j \neq i}^3 u_{j,\kappa} \quad & \mathrm{in} \ \Omega  \\
u_{i,\kappa} \geq 0 & \mathrm{in} \ \Omega \\
u_{i,\kappa}=\gamma_i & \mathrm{on} \ \partial \Omega.\end{array} \right.
\end{equation}
Here $k$ represents the competition between two different densities; the authors analyze the behavior of the solutions as $\kappa \to +\infty$. They prove convergence to a limiting configuration $(u_1,u_2,u_3)$ and presence of the segregation phenomenon, that is $u_i\cdot u_j=0$ a.e. $x\in\Omega$, for $i\neq j$. Our Theorem \ref{theorem_uniqueness_Sgamma} implies that the limiting configuration is in fact solution of (\ref{problema_minimizzazione}), more precisely
\begin{teo} \label{theorem_uniqueness}
Let $(\gamma_i) \in g$ and $V$ satisfy (\ref{hypothesis_V}). Every solution $(u_{1,\kappa}, u_{2,\kappa}, u_{3,\kappa})$ of \eqref{sistema_popolazioni} satisfies
\begin{itemize}
\item[(i)] the whole sequence $u_{i, \kappa}$ converges to a function $u_i$ in $H^1\cap C^{0,\alpha}(\Omega)$ for every $\alpha \in (0,1)$, as $\kappa \to +\infty$;
\item[(ii)] the limiting triple $(u_1,u_2,u_3) \in \Sgamma$ achieves the minimum in (\ref{problema_minimizzazione}) and correspondingly the supports are solutions of (\ref{minimum_pb}).
\end{itemize}
\end{teo}

\section{Continuous dependence of the nodal arcs with respect to the boundary trace}\label{section_continuous_dependence}
Unlike the previous sections, here we let the boundary trace vary, hence it will be convenient to adopt a different notation.

\begin{defi}
We denote by $\Gammaa$ a trace belonging to $\mathcal{G}$, having a triple point at $a$. Let $\Ua$ be the solution of (\ref{cauchy_pb_A}) with boundary trace $\Gammaa$, singularity at $a$ and magnetic potential defined in Lemma \ref{lemma_potenziale2}. Then the nodal set of $\Ua$ consists of three arcs meeting at $a$; we will denote by $\etaa(t):(T_1,T_2)\to \C$ a parametrization of one nodal arc of $\Ua$.
\end{defi}

It is worth to stress that every function $\Ua$ considered in this section has a triple point. Notice that, thanks to Theorem \ref{theorem_local_uniqueness}, all the boundary traces sufficiently close to $\Gammaa$ in the $L^\infty$--norm, also admit triple point. It holds

\begin{teo} \label{regularity_nodalset}
Assume that (\ref{hypothesis_V}) holds. Let $\Omegatilde \subset\subset D$ and let $\beta \in (0,1/2)$. Given $\Gamma_{a_1}$, there exist $\epsilon, C>0$ such that for every $\Gamma_{a_2}\in \mathcal G$ with $\|\Gamma_{a_1}-\Gamma_{a_2}\|_{L^\infty}(\partial D)<\epsilon$, it holds
$$
||\eta_{a_1} -\eta_{a_2}||_{C^{1,\beta}(\Omegatilde)} \leq C ||\Gamma_{a_1}-\Gamma_{a_2}||_{L^\infty(\partial D)},
$$
for a suitable choice of the nodal arcs and of the parametrization.
\end{teo}

The rest of the paragraph is devoted to the proof of Theorem \ref{regularity_nodalset}, hence we will tacitly assume the hypothesis and notations stated therein; in particular $\Omegatilde \subset\subset D$ is fixed. Notice that, by applying the conformal transformation $T_{a_1}$ defined in (\ref{eq:moebius}), we can assume without loss of generality that $a_1=0$, $a_2=a$. As usual, it will be convenient to work with the functions $\uunoa$ introduced in Lemma \ref{lemma_definition_uuno}; we start with some estimates.

\begin{lemma}\label{lemma_C^1_estimates_bis}
Let $\alpha \in (0,1)$. Given $\GammaO$ there exist $\epsilon,C>0$ such that for every $\Gammaa$ satisfying $||\Gammaa-\GammaO||_{L^\infty(\partial D)}<\epsilon$, it holds
\[
\|\Vunoa -\VunoO \|_{L^\infty(\Omegatilde)}+\|\uunoa-\uunoO\|_{C^{1, \alpha}(\Omegatilde)}\leq C||\Gammaa-\GammaO||_{L^\infty(\partial D)}.
\]
\end{lemma}
\begin{proof}
By Theorem \ref{theorem_local_uniqueness} we infer the existence of $\epsilon_1,C_1$ such that $|a|\leq C_1 ||\Gammaa-\GammaO||_{L^\infty(\partial D)}$, for every $\Gammaa$ satisfying $||\Gammaa-\GammaO||_{L^\infty(\partial D)}<\epsilon_1$. Proceeding as in Lemma \ref{lemma_C^1_estimates} we obtain
\[
\|\Vduea -\VdueO \|_{L^\infty(D)}+\|\uduea-\udueO\|_{C^{1, \alpha}(\Omegatilde)}\leq 
C_2 ||\Gammaa-\GammaO||_{L^\infty(\partial D)}.
\]
The same estimate clearly holds for $\uunoa$, provided we choose $\epsilon\leq\epsilon_1$ in such a way that $\Omegatilde\subset \Omegaa$.
\end{proof}

Motiveted by the previous lemma, we define
\[
\Gzero=\{ \Gammaa \in \mathcal G : ||\Gammaa-\GammaO||_{L^\infty(\partial D)}<\epsilon\}.
\]
Therefore we can rewrite the statement of Theorem \ref{regularity_nodalset} in the following way: given $\alpha \in (0,1)$ there exists $C>0$ such that
$$
||\etaunoa -\etaunoO||_{C^{1,\alpha}(\Omegatilde)} \leq C ||\Gammaa-\GammaO||_{L^\infty(\partial D)}, \qquad \forall \ \Gammaa \in \Gzero,
$$
for a suitable choice of the nodal arcs and of the parametrization. Here $\etaunoa$ is any branch of nodal line of $\uunoa$, and by definition it holds $\etaa(t)=(\etaunoa(t))^2$.

The main tool for our analysis will be Theorem \ref{theorem_nodal_set_f}, in particular the H\"older regularity results proved by Helffer, Hoffmann--Ostenhof and Terracini in \cite{HHOT}. For clarity of exposition, we shall recall in the following lemma how these results apply to the function $\uunoa$.

\begin{lemma}\label{lemma_HHOT}
There exists $\xitildea\in C^{0,\a}(\Omegaa,\C), \forall \alpha \in (0,1)$, with $\xitildea(0)=0$, such that
\begin{equation*}
\uunoa (\rho,\phi)=\frac{\rho^{3}}{3}\Big\{ c_{3}(a) \cos(3\phi) + d_{3}(a) \sin(3\phi) +\xitildea(\rho, \phi) \Big\},
\end{equation*}
where $y=\rho e^{i\phi}$. Equivalently, there exists $\xia \in C^{0,\a}(\Omegaa,\C), \forall \alpha \in (0,1)$, such that
\begin{equation}\label{xi_a}
2\parder{\uunoa}{y}(y)= y^2\xia(y), \qquad \xia(0)=c_{3}(a)-i d_{3}(a) \neq 0.
\end{equation}
Moreover for every $k \leq 2$ the following Cauchy formula is available
\begin{equation}\label{gradiente_uduea}
\frac{2}{y^k}\parder{\uunoa}{y}(y)=-\frac{i}{\pi}\int_{\partial \Omegaa} \frac{1}{z^k(z-y)}\parder{\uunoa}{z}(z) \ dz
+\frac{1}{2\pi}\int_{\Omegaa} \frac{-\Delta \uunoa(z)}{z^k(z-y)} \ dz_1dz_2
\end{equation}
where the first integral is a complex line integral, whereas the second one is a double integral in the real variables $z_1,z_2$.
\end{lemma}

Using some ideas in \cite{HHOT}, Theorem 2.1, we can prove the following estimate.

\begin{lemma} \label{lemma_xi_estimates}
Let $\xia, \xitildea$ be the functions defined in the previous lemma. Given $\alpha \in (0,1)$, there exists a positive constant $C$ such that for every $\Gammaa\in \Gzero$ it holds
$$
\|\xia-\xiO\|_{C^{0,\alpha}(\Omegatilde)}\leq C ||\Gammaa-\GammaO||_{L^\infty(\partial D)}.
$$
The same estimate holds for the function $\xitildea$.
\end{lemma}

\begin{proof}
We prove this result by induction on $k$, where for $k=0,1,2$ we define
\begin{eqnarray} 
T_k: & D & \to C^{0,\alpha}(\Omegatilde, \C) \nonumber \\
&a& \mapsto \frac{2}{y^k}\parder{\uunoa}{y}. \nonumber
\end{eqnarray}
Suppose first $k=0$. An integral expression for $T_k(a)$ is given by \eqref{gradiente_uduea}, hence it is enough to show the existence of $C>0$ such that
\begin{multline*}
\left|-\frac{i}{\pi}\int_{\partial \Omegatilde} \left(\parder{\uunoa}{z}-\parder{\uunoO}{z}\right) \left(\frac{1}{z-y_1}-\frac{1}{z-y_2}\right) \ dz \right|+ \\
+ \left| \frac{1}{2\pi}\int_{\Omegatilde} \left(-\Delta\uunoa +\Delta\uunoO \right) \left(\frac{1}{z-y_1}-\frac{1}{z-y_2}\right) 
\ dz_1dz_2 \right| \leq C ||\Gammaa-\GammaO||_{L^\infty(\partial D)} |y_1-y_2|^\alpha,
\end{multline*}
for every $\Gammaa \in \Gzero$. The first integral is smooth in $y$, hence it is sufficient to apply Lemma \ref{lemma_C^1_estimates_bis}. As it concerns the second integral, following \cite{HHOT} we write
\begin{equation*}
\begin{split}
&\left|\int_{\Omegatilde} \left(\Vunoa\uunoa -\VunoO\uunoO \right) \left(\frac{1}{z-y_1}-\frac{1}{z-y_2}\right) \ dz_1dz_2 \right| \\
&\leq \int_{\Omegatilde} \left(|\Vunoa||\uunoa -\uunoO| +|\Vunoa-\VunoO||\uunoO| \right) \left|\frac{1}{z-y_1}-\frac{1}{z-y_2}\right| \ dz_1dz_2 \\
&\leq \left(||\Vunoa||_{L^\infty(\Omegatilde)}||\uunoa -\uunoO||_{L^\infty(\Omegatilde)} +||\Vunoa-\VunoO||_{L^\infty(\Omegatilde)}||\uunoO||_{L^\infty(\Omegatilde)} \right)
\int_{\Omegatilde} \frac{|y_1-y_2|}{|z-y_1||z-y_2|}\ dz_1dz_2 \\
&\leq C||\Gammaa-\GammaO||_{L^\infty(\partial D)} |y_1-y_2| \log|y_1-y_2|,
\end{split}
\end{equation*}
where we used Lemma \ref{lemma_C^1_estimates_bis} in the last inequality. This concludes the proof for $k=0$; now assuming that the result holds for $k=0$ or $k=1$, let us prove it for $k+1$. By assumption there exists $C>0$ such that for every $\Gammaa\in \Gzero$ it holds
$$
\left\|\frac{2}{y^k}\parder{\uunoa}{y}-\frac{2}{y^k}\parder{\uunoO}{y}\right\|_{C^{0,\alpha}(\Omegatilde)} \leq C||\Gammaa-\GammaO||_{L^\infty(\partial D)}.
$$
Since $k<2$ and the origin is a zero of order two for $\uunoa$, we have $\lim_{y \to 0} \frac{1}{y^k}\parder{\uunoa}{y}=0, \ \forall a$. As a consequence, the inductive assumption gives
$$
\sup_{y\in \Omegatilde} \frac{1}{|y|^{k+\alpha}}\left|\parder{\uunoa}{y}(y)-\parder{\uunoO}{y}(y)\right| \leq C ||\Gammaa-\GammaO||_{L^\infty(\partial D)}.
$$
Following \cite{HW} we use the identity
$$
\uunoa(\rho,\phi)= \int^\rho_0 \left( \parder{\uunoa}{y_1}(t,\phi)\cos\phi+\parder{\uunoa}{y_2}(t,\phi)\sin\phi \right) dt
$$
which implies, together with the previous inequality
\begin{equation}\label{holder1}
|(\uunoa-\uunoO)(y)|\leq \int^1_0 |y| \left|2\left( \parder{\uunoa}{y}-\parder{\uunoO}{y} \right)(ty) \right| dt \leq  C||\Gammaa-\GammaO||_{L^\infty(\partial D)}|y|^{k+1+\alpha},
\end{equation}
for every $y\in \Omegatilde$. Now we can proceed as in the case $k=0$:
\begin{equation*}
\begin{split}
& \left| \int_{\Omegatilde} \frac{-\Delta\uunoa +\Delta\uunoO}{z^{k+1}} \left(\frac{1}{z-y_1}-\frac{1}{z-y_2}\right) 
\ dz_1dz_2 \right| \\
&\leq \left(||\Vunoa||_{L^\infty(\Omegatilde)} \left\| \frac{\uunoa -\uunoO}{z^{k+1}}\right\|_{L^\infty(\Omegatilde)} +||\Vunoa-\VunoO||_{L^\infty(\Omegatilde)}\left\|\frac{\uunoO}{z^{k+1}}\right\|_{L^\infty(\Omegatilde)} \right)
|y_1-y_2| \log|y_1-y_2|, \\
&\leq C ||\Gammaa-\GammaO||_{L^\infty(\partial D)}|y_1-y_2| \log|y_1-y_2|,
\end{split}
\end{equation*}
where we used \eqref{holder1} in the last inequality. This conclude the estimate for $\xia$, it is easy to see that it implies the same estimate on $\xitildea$.
\end{proof}

Choose now a branch of nodal line of $\uunoa$ satisfying
\[
\etaunoa:(T_1,T_2) \rightarrow \C, \qquad \lim_{t\to T_1}\etaunoa(t)=0, \qquad \etaunoa(t) \subset \Omegatilde \quad \forall \ t\in(T_1,T_2),
\]
with $T_1,T_2$ eventually infinite. Then it must hold
\[
\dot\eta^{(2)}_{a}(t)=- i \kappa_a(t) \parder{\uunoa}{\bar{y}}(\etaunoa(t)), \qquad \uunoa(\etaunoa(t))=0,
\]
where $\kappa_a(t)$ is any real function, sufficiently regular in $(T_1,T_2)$. Since every $\uunoa$ has a zero of order two at the origin, we choose
\[
\kappa_a(t)=\frac{1}{|\etaunoa(t)|^2}.
\]
By using the polar coordinates $\etaunoa(t)=\rhoa(t) e^{i \phia(t)}$, the equation for the curve becomes
\begin{eqnarray}\label{eq:nodal_lines_polar}
\left\{ \begin{array}{ll}
\dot{\rhoa}=\frac{1}{|y|^{3}}\left( y_1 \parder{\uunoa}{y_2}- y_2 \parder{\uunoa}{y_1}  \right)
= - \frac{1}{|y|^{3}} \Im(y^3 \xia(y)) \\
\dot{\phia}=-\frac{1}{|y|^{4}}\left( y_1 \parder{\uunoa}{y_1}+ y_2 \parder{\uunoa}{y_2}  \right)
=-\frac{1}{|y|^{4}} \Re(y^3 \xia(y)).
\end{array} \right.
\end{eqnarray}

\begin{lemma}\label{lemma_nodal_line_velocity}
With this choice of the parametrization the interval $(T_1,T_2)$ is bounded, in particular we can choose $T_1=0$.
\end{lemma}

\begin{proof}
By computing the velocity of the curve we obtain
\begin{equation*}
|\dot{\eta}^{(1)}_a(t)|=\frac{1}{|\etaunoa(t)|^2}\left|\parder{\uunoa}{\bar{y}}(\etaunoa(t))\right|
=|\xia(\etaunoa(t))| \to |\xia(0)| =\sqrt{c_{3}(a)^2+d_{3}(a)^2} \ 
\qquad \text{as } t \to T_1,
\end{equation*}
where $c_3(a),d_3(a)$ are the first nontrivial terms in the asymptotic expansion (\ref{xi_a}). Hence $t$ is asymptotically a multiple of the arc length as $t \to T_1$ and the lemma is proved.
\end{proof}

\begin{proof}[End of the proof of Theorem \ref{regularity_nodalset}]
By writing the equation of the curve in polar coordinates, the time derivative writes $\dot{\eta}^{(1)}_a(t)=e^{i\phia(t)}(\dot{\rhoa}(t)+i\rhoa(t)\dot{\phia}(t))$. Therefore we wish to show that, for every $\alpha\in(0,1)$ there exist constants $K_1,K_2$ such that
\[
\begin{split}
&\| \dot \rhoa - \dot \rhoO \|_{C^{0,\alpha}([0,T_2])} +
\| \rhoa\dot\phia - \rhoO\dot\phiO \|_{C^{0,\alpha}([0,T_2])}
\leq K_1 ||\Gammaa-\GammaO||_{L^\infty(\partial D)},\\
&\| \phia-\phiO \|_{C^{0,\alpha}([0,T_2])} \leq K_2 ||\Gammaa-\GammaO||_{L^\infty(\partial D)},
\end{split}
\]
for every $\Gammaa \in \Gzero$.

The first inequality comes directly from equation (\ref{eq:nodal_lines_polar}) and Lemma \ref{lemma_xi_estimates}, by regularity results for ordinary differential equations. Let us prove the second one; from $\uunoa(\etaunoa(t))=0$ we deduce
$$
c_3(a) \cos(3\phia) + d_{3}(a) \sin(3\phia) +\xitildea(\etaa(t))=0,
$$
and hence
\begin{eqnarray*}
\dot{\phia} &=&-\frac{1}{|y|} (c_{3}(a) \cos(3\phi) + d_{3}(a) \sin(3\phi) )
-\frac{1}{|y|^4}\Re[y^3 (\xia(y)-\xia(0))]\\
&=&\frac{\xitildea(y)}{|y|}-\frac{1}{|y|^4}\Re[y^3 (\xia(y)-\xia(0))].
\end{eqnarray*}
Now, since both $\xitildea(y)$ and $\xia(y)-\xia(0)$ satisfy Lemma \ref{lemma_xi_estimates} and vanish at the origin, there exists $C>0$ such that
\begin{equation}\label{zeta_xitilde_holder}
\max \left\{ \sup_{y \in \Omegatilde} \frac{|\xia(y)-\xia(0)|}{|y|^\alpha}, 
\sup_{y \in \Omegatilde} \frac{|\xitildea(y)|}{|y|^\alpha} \right\} \leq C ||\Gammaa-\GammaO||_{L^\infty(\partial D)},
\end{equation}
for every $\Gammaa \in \Gzero$. By combining these results with Lemma \ref{lemma_nodal_line_velocity}, we obtain
\[
|\dot{\phia}(t)-\dot{\phiO}(t)| \leq C ||\Gammaa-\GammaO||_{L^\infty(\partial D)} t^{\alpha-1},
\]
and finally $\|\phia-\phiO\|_{C^{0,\alpha}([0,T_2])} \leq K_2 ||\Gammaa-\GammaO||_{L^\infty(\partial D)}$, which concludes the proof of the theorem.
\end{proof}


\noindent Universit\`a di Milano Bicocca,\\
Dipartimento di Matematica e Applicazioni,\\
Via R. Cozzi 53, 20125 Milano, Italy.\\
\textit{E-mail addresses:} \texttt{b.noris@campus.unimib.it, susanna.terracini@unimib.it}

\end{document}